\documentclass[11pt]{amsart}
\usepackage{amsmath,amssymb,mathtools}
\usepackage{microtype}
\usepackage[a4paper,margin=29mm,headheight=14pt]{geometry}
\usepackage[hidelinks,pdfauthor={Chang Li}]{hyperref}
\numberwithin{equation}{section}
\newtheorem{theorem}{Theorem}[section]
\newtheorem{lemma}[theorem]{Lemma}
\newtheorem{proposition}[theorem]{Proposition}
\newtheorem{corollary}[theorem]{Corollary}
\theoremstyle{definition}
\newtheorem{definition}[theorem]{Definition}
\theoremstyle{remark}
\newtheorem{remark}[theorem]{Remark}
\allowdisplaybreaks[1]
\date{}
\title{Carleson duality and endpoint structure of Euclidean flag Triebel--Lizorkin spaces}
\author[C. Li]{Chang Li}
\address{Department of Mathematics, Shanghai University,
Shanghai 200444, People's Republic of China}
\email{changli@shu.edu.cn}

\begin{document}
\begin{abstract}
The main purpose of this paper is to develop the duality theory of
the homogeneous Triebel--Lizorkin spaces associated with the
Euclidean flag, using the discrete Littlewood--Paley--Stein
analysis. After establishing bounded analysis and synthesis maps, we
describe the continuous duals by Carleson conditions. In the Hardy
spatial range, the representing distributions satisfy an open-set
Carleson condition when the sequence exponent is greater than one and
a coefficient supremum condition when the sequence exponent is at
most one; in the Banach spatial range the dual is a
conjugate-exponent or a supremum space. At the supremum sequence
endpoint, the closure of the flag test functions is characterized by
a pointwise scale tail measured in $L^p$, and a coefficient pairing
decomposes every functional on the full space into a regular part
and a part that annihilates the closure. By means of a
fixed-envelope modulation construction we show that this closure is
proper for every finite spatial exponent and that the annihilator is
nonzero when the spatial exponent is at least one. The same
analysis and synthesis framework also yields complex interpolation,
real interpolation in the spatial exponent, the lifting property and
the standard mixed-norm embeddings.
\end{abstract}
\maketitle

\section{Introduction and statement of results}
\label{sec:introduction}

The duality of function spaces adapted to singular integrals depends
on the underlying dilations and on the order of the spatial and
scale norms. In the one-parameter setting, the Besov and
Triebel--Lizorkin scales provide a systematic description of
smoothness and integrability; see Triebel \cite{Triebel83}. In the
multiparameter setting, product and flag geometries lead to different
cancellation conditions and endpoint spaces.

The product theory is governed by independent dilations in the
coordinate groups. Journ\'e \cite{Journe85} developed the $T1$ theory
for non-convolution product operators, and Chang and Fefferman
\cite{CF80} identified the dual of the product Hardy space by an
open-set Carleson condition. These results extend aspects of the
one-parameter Hardy-space theory of Fefferman and Stein
\cite{FS72}. Discrete transforms of the type introduced by Frazier
and Jawerth \cite{FJ85,FJ90} provide the sequence-space framework for
studying Besov and Triebel--Lizorkin norms.

Flag singular integrals arise from the Marcinkiewicz multiplier
theory on the Heisenberg group developed by M\"uller, Ricci and Stein
\cite{MRS95,MRS96}. The relevant multipliers have a two-parameter
scaling structure, although the corresponding two-parameter
dilations do not act by automorphisms of the Heisenberg group. The canonical
Euclidean model is the two-step flag on $\mathbb R^n\times\mathbb R^m$,
\begin{equation}\label{eq:flag}
 \{0\}\subset\{0\}\times\mathbb R^m\subset\mathbb R^n\times\mathbb R^m.
\end{equation}
Its kernels were written down by Nagel, Ricci and Stein
\cite{NRS01}; their singularity is concentrated on the set
$\{0\}\times\mathbb R^m$, a hybrid of the point and product
singularities, and the construction was carried to homogeneous
groups in \cite{NRSW12}.

The endpoint theory of flag singular integrals requires spaces
adapted to their partial cancellation: the classical and product
Hardy spaces do not capture the flag geometry. The appropriate
intermediate space, the flag Hardy space, was developed through the
discrete Littlewood--Paley--Stein analysis of Han and Lu
\cite{HL08,HL10}, with later real-variable forms in
\cite{HLLW22,HLW19}; the underlying reproducing-formula techniques
go back to Han and Sawyer \cite{HS94} and Han \cite{Han00}. On the
function-space side, Ding, Lu and Ma \cite{DLM10} introduced the
multiparameter flag Besov and Triebel--Lizorkin scales, the scales
on which we work throughout.

Endpoint duality in a multiparameter setting is naturally described
by Carleson conditions. In the product case, the Chang--Fefferman
condition is taken over open sets, and the covering arguments of
Journ\'e \cite{Journe85} relate these conditions to the geometry of
rectangles. Duality theorems are also available for related flag
structures: Yang \cite{Yang09} considered a parabolic geometry, and
Ding and Lu \cite{DL16} treated an isotropic--parabolic composition,
including quasi-Banach and finite-Banach ranges. Here we consider
the Euclidean flag in \eqref{eq:flag}, whose first decomposition is
isotropic and whose second decomposition acts only in $y$. We give
a unified dual classification for arbitrary real smoothness and
finite positive spatial and sequence exponents, and treat the
supremum sequence endpoint separately. The case $s=0$, $q=2$ and
$0<p\leq1$ is already contained in the flag Hardy theory of Han, Lin
and Wu \cite{HLW19}. At $q=\infty$, we identify the dual of the
closure of the test functions by coefficient pairings; on the full
space, these pairings give the regular part of each continuous
functional, with the remaining part annihilating that closure.

The discrete analysis--synthesis framework for the Euclidean flag
is shared with the companion paper \cite{LiBesov}, which develops
coefficient block duality and interpolation for the flag Besov
spaces. That paper is available as
\href{https://arxiv.org/abs/2610.32204v1}{arXiv:2610.32204v1};
all references to its result numbers refer to this version.
The test-space identification and exact sampling formula are
recalled in the form needed here, with precise references to
\cite[Proposition~1.1 and Section~2]{LiBesov}. We include proofs of
the pointwise grid, Peetre and supremum sampling estimates shared
with \cite[Lemmas~2.2, 2.4 and~2.5]{LiBesov}, so that the estimates
used in the Triebel--Lizorkin argument can be checked directly.
The vector-valued maximal estimates, the open-set Carleson duality,
the fixed-envelope modulation construction and the spatial-index
interpolation are established in the present paper for the
$L^p(\ell^q)$ ordering.

In this paper we give a unified coefficient classification for the
Euclidean flag, and we give a careful treatment of the supremum
sequence endpoint. In particular, we distinguish the closure of the
test functions from the full $q=\infty$ space, identify the dual of
that closure by means of absolutely convergent coefficient pairings,
and, for $1\leq p<\infty$, exhibit nonzero functionals on the full
space that vanish on the closure. The fixed-envelope modulation construction in
Proposition~\ref{prop:modulations} shows that the inclusion is strict
in the Triebel--Lizorkin norm. Maximal-function and difference
characterisations of inhomogeneous flag spaces were obtained by Wu
and Liu \cite{WL12}; relations with product spaces were studied by
Chang, Han and Wu \cite{CHW21} and by Cao, Chang and Wu
\cite{CCW24}; and a related three-parameter setting was considered
by Chen, Tao and Zheng \cite{CTZ25}.

\subsection{Spaces and dyadic frames}

A flag rectangle is a set $R=I\times J$, where $I\subset\mathbb R^n$
and $J\subset\mathbb R^m$ are cubes with sides parallel to the axes
and $\ell(I)\leq\ell(J)$. The flag maximal function is
\begin{equation}\label{eq:MF}
 M_F h(z)=\sup_{\substack{R\ni z\\R\text{ a flag rectangle}}}
                         \frac1{|R|}\int_R|h(u)|\,du.
\end{equation}
\begin{definition}[\cite{DLM10}]
Call a Schwartz function $F(x,y,z)$ on
$\mathbb R^n\times\mathbb R^m\times\mathbb R^m$ a product test
function on $\mathbb R^{n+m}\times\mathbb R^m$ when its moments
vanish in each of the two test groups, that is,
\begin{equation*}
 \int_{\mathbb R^{n+m}}F(x,y,z)x^\alpha y^\beta\,dx\,dy=0,
 \qquad
 \int_{\mathbb R^m}F(x,y,z)z^\gamma\,dz=0,
\end{equation*}
for every choice of multi-indices $\alpha,\beta,\gamma$, pointwise in
the remaining variables.
\end{definition}

Write $\mathcal S_\infty(\mathbb R^{n+m}\times\mathbb R^m)$ for the
space of product test functions, equipped with the topology of the
Schwartz seminorms. The flag test space is the projection of this
space under the map that integrates out the last group of variables.

\begin{definition}[\cite{DLM10}]
A function $f$ on $\mathbb R^n\times\mathbb R^m$ belongs to the flag
test space $\mathcal S_F$ if there is a product test function
$F\in\mathcal S_\infty(\mathbb R^{n+m}\times\mathbb R^m)$ such that
\begin{equation}\label{eq:piF}
 f(x,y)=\int_{\mathbb R^m}F(x,y-z,z)\,dz.
\end{equation}
Give $\mathcal S_F$ the quotient locally convex topology carried by
\eqref{eq:piF}: a net converges in $\mathcal S_F$ exactly when it has
lifts converging in $\mathcal S_\infty$. In terms of the Schwartz
seminorms, the same topology is generated by the induced seminorms
\[
 \rho_F(f)=\inf\{\rho(F):F\in\mathcal S_\infty
 \text{ is a lift of }f\text{ through }\eqref{eq:piF}\}.
\]
\end{definition}

The continuous dual is denoted by $\mathcal S_F'$. Pairings are
linear in the distribution and conjugate linear in the test
function. For later use we record the partial-cancellation
description of this space, the one adopted in
\cite[Definition~1.2]{HLW19}:
\begin{equation}\label{eq:partialmoments}
 \mathcal S_F=\Bigl\{f\in\mathcal S(\mathbb R^{n+m}):
 \int_{\mathbb R^m}f(x,y)y^\beta\,dy=0
 \ \text{for every }x\in\mathbb R^n
 \text{ and every multi-index }\beta\Bigr\},
\end{equation}
and the quotient topology agrees with the subspace topology inherited
from $\mathcal S(\mathbb R^{n+m})$. The equivalence follows from the
Fourier-side identity $\widehat{\pi F}(\xi,\eta)=\widehat F(\xi,\eta,\eta)$
together with a flat right inverse on the product side; the cutoff
and seminorm estimates are written out in the companion paper
\cite[Proposition~1.1]{LiBesov} and are not repeated here.

Hahn--Banach then extends every functional on this closed subspace
to a tempered distribution, so $\mathcal S_F'$ identifies canonically
with the quotient of $\mathcal S'$ by the annihilator of
$\mathcal S_F$. In particular, convolving with a translate of an
annular flag analysing function does not depend on the chosen
tempered representative, because every such translate lies in
$\mathcal S_F$.

Choose real even Schwartz functions $\psi^{(1)}\in\mathcal S(\mathbb R^{n+m})$
and $\psi^{(2)}\in\mathcal S(\mathbb R^m)$, with Fourier transforms
supported in the annulus $\{1/2\leq|\cdot|\leq2\}$ in their
respective spaces, whose moments of every order therefore vanish.
Suppose that
\begin{align}
 \sum_{j\in\mathbb Z}
  |\widehat{\psi^{(1)}}(2^{-j}\xi,2^{-j}\eta)|^2&=1
       &&((\xi,\eta)\ne0),\label{eq:partition1}\\
 \sum_{k\in\mathbb Z}
  |\widehat{\psi^{(2)}}(2^{-k}\eta)|^2&=1
       &&(\eta\ne0).\label{eq:partition2}
\end{align}
Such a pair is called admissible.  Put
\begin{equation}\label{eq:psijk}
 \psi_{j,k}(x,y)
 =\psi^{(1)}_{j}\ast_2\psi_{k}^{(2)}=\int_{\mathbb R^m}2^{j(n+m)}\psi^{(1)}
        \bigl(2^jx,2^j(y-z)\bigr)\,
                  2^{km}\psi^{(2)}(2^kz)\,dz .
\end{equation}
Here $\ast_2$ denotes convolution in the second variable. We refer
to $\psi_{j,k}$ as the flag analysing kernels. The notation
$\varphi_{j,k}$ has the same meaning for another admissible pair
$\varphi^{(1)},\varphi^{(2)}$. Their Fourier supports satisfy
\begin{equation}\label{eq:flag-support-explicit}
 \operatorname{supp}\widehat\psi_{j,k}
 \subset\bigl\{2^{j-1}\leq|(\xi,\eta)|\leq2^{j+1},\quad
                2^{k-1}\leq|\eta|\leq2^{k+1}\bigr\}.
\end{equation}
Consequently $\psi_{j,k}=0$ unless $k\leq j+1$; at the possible
boundary $k=j+2$, smoothness and annular support force the product
multiplier to vanish.

For $s=(s_1,s_2)\in\mathbb R^2$, we use the smoothness weight
$2^{-js_1}2^{-(j+k)s_2}=2^{-j(s_1+s_2)}2^{-ks_2}$.
On the support of a nonzero kernel,
$|(\xi,\eta)|\asymp2^j$ and $|\eta|\asymp2^k$; thus the smoothness
parameters attached to these two frequency scales are $s_1+s_2$
and $s_2$, respectively, with the negative-sign convention in the
norm below. The pair $s$ is a linear reparametrisation of those
scale parameters.

\begin{definition}\label{def:F}
Let $s\in\mathbb R^2$, $0<p<\infty$ and $0<q\leq\infty$. The
homogeneous flag Triebel--Lizorkin space
$\dot F^{s,q}_{p,F}=\dot F^{s,q}_{p,F}(\mathbb R^n\times\mathbb R^m)$
consists of those $f\in\mathcal S_F'$ for which
\begin{equation}\label{eq:Fdefinition}
 \|f\|_{\dot F^{s,q}_{p,F}}
 =\left\|\left\{\sum_{j,k\in\mathbb Z}
 \bigg(2^{-js_1}2^{-(j+k)s_2}|\psi_{j,k}*f|\bigg)^q
                  \right\}^{1/q}\right\|_{L^p}<\infty,
\end{equation}
with the inner expression a supremum when $q=\infty$.
\end{definition}

The effective lengths and the sampling grids are
\begin{equation}\label{eq:grid}
 a_j=2^{-j},\qquad b_{j,k}=2^{-j}\vee2^{-k},\qquad
 \mathcal R_{j,k}^N
 =\{I\times J:\ell(I)=2^{-N}a_j,\
                    \ell(J)=2^{-N}b_{j,k}\},
\end{equation}
where $I,J$ range over the usual half-open dyadic cubes and $z_R$
denotes the lower corner of $R$. Throughout, $N$ is a fixed integer
with $N\geq N_0\geq0$, where $N_0$ is chosen sufficiently large for the
sampling estimates below, including the enlarged support box used
in regularisation. For any finite collection of analysing systems
or endpoint parameters, we use one common such $N$. At a fixed
scale each rectangle has volume
\begin{equation}\label{eq:Rvolume}
 |R|=2^{-N(n+m)}2^{-jn}2^{-m\min(j,k)}.
\end{equation}
For the active indices $k\leq j+1$, this gives
$|R|\asymp_N 2^{-jn}2^{-km}$; this comparison is used only on
those indices, since $\psi_{j,k}=0$ when $k>j+1$.

\begin{definition}\label{def:f}
For a labelled sequence $t=\{t_{j,k,R}\}$ the coefficient space
$\dot f^{s,q}_{p,F}$ is defined by
\begin{equation}\label{eq:fsequence}
 \|t\|_{\dot f^{s,q}_{p,F}}
 =\left\|\left\{\sum_{j,k}\sum_{R\in\mathcal R_{j,k}^N}
       \bigg(2^{-js_1}2^{-(j+k)s_2}|t_{j,k,R}|\bigg)^q\mathbf1_R
 \right\}^{1/q}\right\|_{L^p},
\end{equation}
with the inner expression a supremum when $q=\infty$.
\end{definition}

We use unnormalised coefficients, so that the factor $|R|$ occurs in
synthesis.

\begin{definition}\label{def:ST}
The analysis and synthesis maps associated with an admissible pair
are
\begin{equation}\label{eq:analysis-synthesis}
 S_\psi f=\{(\psi_{j,k}*f)(z_R)\}_{j,k,R},\qquad
 T_\psi t=\sum_{j,k}\sum_{R\in\mathcal R_{j,k}^N}
       |R|\,t_{j,k,R}\psi_{j,k}(\,\cdot-z_R).
\end{equation}
\end{definition}

\begin{theorem}[Retraction and discrete norm]\label{thm:retraction}
Let $s\in\mathbb R^2$, $0<p<\infty$ and $0<q\leq\infty$. Then
\[
 S_\psi:\dot F^{s,q}_{p,F}\longrightarrow\dot f^{s,q}_{p,F},
 \qquad
 T_\psi:\dot f^{s,q}_{p,F}\longrightarrow\dot F^{s,q}_{p,F}
\]
are bounded, $T_\psi S_\psi=I$, and
\begin{equation}\label{eq:discrete-norm}
 \|f\|_{\dot F^{s,q}_{p,F}}
 \asymp\|S_\psi f\|_{\dot f^{s,q}_{p,F}}
 \asymp\inf_{T_\psi t=f}\|t\|_{\dot f^{s,q}_{p,F}}.
\end{equation}
Every synthesis series converges weakly in $\mathcal S_F'$; if
$q<\infty$ it also converges in $\dot F^{s,q}_{p,F}$. The space
$\dot F^{s,q}_{p,F}$ is quasi-Banach, continuously embedded in
$\mathcal S_F'$, and is Banach when $p,q\geq1$. The test space
$\mathcal S_F$ is continuously contained in it and is dense when
$q<\infty$.
\end{theorem}

\subsection{The little space and the supremum endpoint}

Write $J=(j,k)$ and $|J|_1=|j|+|k|$. For
$0<q\leq\infty$ define the scale tail
\begin{equation}\label{eq:omega}
 \omega^F_{f,q}(K)
 =\bigl\|\{2^{-js_1}2^{-(j+k)s_2}|\psi_{j,k}*f|\}_{|J|_1>K}\bigr\|_{L^p(\ell^q)},
\end{equation}
with the inner expression a supremum when $q=\infty$, and put
$d_F(f)=\lim_{K\to\infty}\omega^F_{f,\infty}(K)$.

\begin{theorem}[The little Triebel--Lizorkin space]\label{thm:little}
Let $s\in\mathbb R^2$ and $0<p<\infty$. For
$f\in\dot F^{s,\infty}_{p,F}$,
\begin{equation}\label{eq:distance}
 d_F(f)\leq
 \inf_{g\in\mathcal S_F}\|f-g\|_{\dot F^{s,\infty}_{p,F}}
 \leq C d_F(f).
\end{equation}
Consequently the little space
$\dot F^{s,\infty}_{p,F,0}:=\overline{\mathcal S_F}^{\,\dot F^{s,\infty}_{p,F}}$
consists precisely of those $f$ satisfying
\begin{equation}\label{eq:little}
 \lim_{K\to\infty}\left\|\sup_{|j|+|k|>K}
       2^{-js_1}2^{-(j+k)s_2}|\psi_{j,k}*f|\right\|_{L^p}=0.
\end{equation}
For members of the ambient space this is equivalent to
\begin{equation}\label{eq:pointwise-little}
 2^{-js_1}2^{-(j+k)s_2}|\psi_{j,k}*f(z)|\longrightarrow0
 \quad\text{as }|j|+|k|\to\infty,
 \quad\text{for almost every }z.
\end{equation}
The little space is the closure of the finite syntheses using the
fixed pair and grid, and its definition is independent of the
admissible pair and the grid.
\end{theorem}

\subsection{Carleson duality}

\begin{definition}\label{def:carleson}
For $\sigma\in\mathbb R^2$ and $1<r<\infty$, let
$\mathrm{CMO}_{F}^{\sigma,p,r}$ consist of those
$g\in\mathcal S_F'$ for which
\begin{equation}\label{eq:carleson-norm}
 \|g\|_{\mathrm{CMO}_{F}^{\sigma,p,r}}
 =\sup_{\Omega}|\Omega|^{1-1/p-1/r}
 \left\{\sum_{j,k}\sum_{\substack{R\in\mathcal R_{j,k}^N\\R\subset\Omega}}
 |R|\bigg(2^{-j\sigma_1}2^{-(j+k)\sigma_2}|(\psi_{j,k}*g)(z_R)|\bigg)^r
 \right\}^{1/r}<\infty,
\end{equation}
where the supremum runs over open sets $\Omega\subset\mathbb R^{n+m}$
with $0<|\Omega|<\infty$. Write $\dot c_F^{\sigma,p,r}$ for the
sequence space obtained by replacing
$(\psi_{j,k}*g)(z_R)$ with $d_{j,k,R}$ in
\eqref{eq:carleson-norm}.
\end{definition}

At $\sigma=0$, $r=2$, \eqref{eq:carleson-norm} is the flag Carleson
condition. The factors $|\Omega|^{1-1/p-1/r}$ and $|R|$ are fixed by
the coefficient pairing and the order of the spatial and sequence
norms.

\begin{theorem}[Carleson duality in the Hardy range]\label{thm:carleson-dual}
Let $s\in\mathbb R^2$, $0<p\leq1$ and $1<q<\infty$. Then
\begin{equation}\label{eq:carleson-dual}
 \bigl(\dot F^{s,q}_{p,F}\bigr)^*
       =\mathrm{CMO}_{F}^{-s,p,q'}
\end{equation}
with equivalent norms, under the unique continuous extension of the
distribution pairing on $\mathcal S_F$. The Carleson space is Banach
and its norm is independent, up to equivalence, of the admissible
pair and the exact sampling grid.
\end{theorem}

\begin{theorem}[Duality in the Banach range]\label{thm:duality}
Let $1<p,q<\infty$. Then
\begin{equation}\label{eq:Fduality}
 \bigl(\dot F^{s,q}_{p,F}\bigr)^*
       =\dot F^{-s,q'}_{p',F}
\end{equation}
with equivalent norms, under the unique continuous extension of the
distribution pairing on $\mathcal S_F$.
\end{theorem}

For $0<p\leq1$ the Carleson notation extends to the endpoints
$r=1$ and $r=\infty$ by
\begin{align}
 \|d\|_{\dot c_F^{\sigma,p,1}}
 &=\sup_{\Omega}|\Omega|^{-1/p}
   \sum_J\sum_{\substack{R\in\mathcal R_J^N\\R\subset\Omega}}
                 |R|2^{-j\sigma_1}2^{-(j+k)\sigma_2}|d_{J,R}|,\label{eq:c-one}\\
 \|d\|_{\dot c_F^{\sigma,p,\infty}}
 &=\sup_{J,R}|R|^{1-1/p}2^{-j\sigma_1}2^{-(j+k)\sigma_2}|d_{J,R}|,\label{eq:c-infinity}
\end{align}
and $\mathrm{CMO}_{F}^{\sigma,p,r}$, for $r\in\{1,\infty\}$, is
defined by $g\in\mathcal S_F'$ and $S_\psi g\in\dot c_F^{\sigma,p,r}$.

\begin{theorem}[Hardy-range endpoint duals]\label{thm:Hardy-endpoint}
Let $0<p\leq1$ and $s\in\mathbb R^2$. For $0<q\leq1$,
\begin{equation}\label{eq:small-q-dual}
 \bigl(\dot F^{s,q}_{p,F}\bigr)^*
             =\mathrm{CMO}_{F}^{-s,p,\infty},
\end{equation}
and
\begin{equation}\label{eq:little-Hardy-dual}
 \bigl(\dot F^{s,\infty}_{p,F,0}\bigr)^*
             =\mathrm{CMO}_{F}^{-s,p,1}.
\end{equation}
The norms are equivalent and the pairing is the unique continuous
extension of the distribution pairing on $\mathcal S_F$. Both
endpoint Carleson spaces are Banach and their norms are independent,
up to equivalence, of the admissible pair and the exact grid.
\end{theorem}

\begin{theorem}[Banach-range endpoint duals]\label{thm:Banach-endpoint}
Let $1<p<\infty$, $s\in\mathbb R^2$ and $0<q\leq1$. Then
\begin{align}
 \bigl(\dot F^{s,q}_{p,F}\bigr)^*
       &=\dot F^{-s,\infty}_{p',F},\label{eq:q-one-dual}\\
 \bigl(\dot F^{s,\infty}_{p,F,0}\bigr)^*
       &=\dot F^{-s,1}_{p',F}.\label{eq:little-Banach-dual}
\end{align}
The pairings are the unique continuous extensions of the
distribution pairing on $\mathcal S_F$ and the norms are equivalent.
\end{theorem}

\begin{corollary}[The $q\leq1$ Hardy dual as a Besov space]\label{cor:small-q-Besov}
For $0<p\leq1$, $0<q\leq1$ and $s\in\mathbb R^2$, put
\begin{equation}\label{eq:Besov-shift}
 \delta=1/p-1,\qquad
 t=\bigl(-s_1+(m-n)\delta,\,-s_2-m\delta\bigr).
\end{equation}
Then
\begin{equation}\label{eq:small-q-Besov}
 \bigl(\dot F^{s,q}_{p,F}\bigr)^*=\dot B^{t,\infty}_{\infty,F},
\end{equation}
where
$\|g\|_{\dot B^{t,\infty}_{\infty,F}}
       =\sup_{j,k}2^{-jt_1}2^{-(j+k)t_2}\|\psi_{j,k}*g\|_{L^\infty}$.
The norms are equivalent and independent of the admissible pair.
\end{corollary}

\subsection{The full supremum dual}

\begin{proposition}[Regular pairing and decomposition]\label{prop:full-supremum}
Put $X=\dot F^{s,\infty}_{p,F}$ and
$X_0=\dot F^{s,\infty}_{p,F,0}$, with $0<p<\infty$, and write
\begin{equation}\label{eq:Y}
 Y=\begin{cases}
  \mathrm{CMO}_{F}^{-s,p,1},&0<p\leq1,\\
  \dot F^{-s,1}_{p',F},&1<p<\infty.
 \end{cases}
\end{equation}
For $f\in X$ and $g\in Y$ the pairing
\begin{equation}\label{eq:regular-pairing}
 \langle f,g\rangle_F=[S_\psi f,S_\psi g]
 =\sum_{J,R}|R|(\psi_J*f)(z_R)\overline{(\psi_J*g)(z_R)}
\end{equation}
converges absolutely, is independent of the admissible pair and the
exact grid, and agrees with the distribution pairing when either
entry is a test function. Every $\ell\in X^*$ has a unique
decomposition
\begin{equation}\label{eq:regular-singular}
 \ell=\ell_g+\ell_0,\qquad
 \ell_g(f)=\langle f,g\rangle_F,\quad g\in Y,\quad
 \ell_0\in X_0^\perp:=\{L\in X^*:L|_{X_0}=0\},
\end{equation}
and
\begin{equation}\label{eq:regular-singular-norm}
 \|\ell\|\asymp\|g\|_Y+\|\ell_0\|.
\end{equation}
\end{proposition}

\begin{proposition}[A proper closure and a nonzero annihilator]\label{prop:modulations}
Let $s\in\mathbb R^2$ and $0<p<\infty$, with
$X=\dot F^{s,\infty}_{p,F}$ and $X_0=\dot F^{s,\infty}_{p,F,0}$.
There is a linear map $E:\ell^\infty(\mathbb N)\to X$ such that
\begin{align}
 \|Eb\|_X&\asymp\|b\|_{\ell^\infty},\label{eq:modulation-norm}\\
 \operatorname{dist}_X(Eb,X_0)
       &\asymp\limsup_{\ell\to\infty}|b_\ell|.\label{eq:modulation-distance}
\end{align}
Hence $Eb\in X_0$ if and only if $b\in c_0$; $X_0$ is a proper closed
subspace of $X$ and $X$ is nonseparable. If $1\leq p<\infty$, then
$X_0^\perp\ne\{0\}$, so in this range some continuous functionals on
$X$ are not represented by the regular pairing of
Proposition~\ref{prop:full-supremum}.
\end{proposition}

\subsection{Interpolation}

\begin{theorem}[Complex interpolation]\label{thm:complex-interpolation}
Let $1<p_0,p_1,q_0,q_1<\infty$, $s^0,s^1\in\mathbb R^2$ and
$0<\theta<1$. Put
\begin{equation}\label{eq:interpolated-parameters}
 \frac1p=\frac{1-\theta}{p_0}+\frac\theta{p_1},\qquad
 \frac1q=\frac{1-\theta}{q_0}+\frac\theta{q_1},\qquad
 s=(1-\theta)s^0+\theta s^1.
\end{equation}
Then
\begin{equation}\label{eq:complex-interpolation}
 \big[\dot F^{s^0,q_0}_{p_0,F},
       \dot F^{s^1,q_1}_{p_1,F}\big]_\theta
                       =\dot F^{s,q}_{p,F}
\end{equation}
with equivalent norms.
\end{theorem}

\begin{theorem}[Real interpolation of the spatial exponent]\label{thm:real-interpolation}
Let $1<p_0<p_1<\infty$, $1\leq q\leq\infty$, $s\in\mathbb R^2$ and
$0<\theta<1$, and set
$1/p=(1-\theta)/p_0+\theta/p_1$. Then
\begin{equation}\label{eq:real-interpolation}
 (\dot F^{s,q}_{p_0,F},\dot F^{s,q}_{p_1,F})_{\theta,p}
                         =\dot F^{s,q}_{p,F}
\end{equation}
with equivalent norms. More precisely, writing
$G_f(u)=\|\{2^{-js_1}2^{-(j+k)s_2}(\psi_J*f)(u)\}_J\|_{\ell^q}$, the continuous
$K$-functional satisfies
\begin{equation}\label{eq:scalar-K}
 K(t,f;\dot F^{s,q}_{p_0,F},\dot F^{s,q}_{p_1,F})
       \asymp K(t,G_f;L^{p_0},L^{p_1}).
\end{equation}
This real interpolation is restricted to the Banach spatial range
and fixed smoothness; the quasi-Banach spatial range and the
interpolation of the smoothness are not addressed.
\end{theorem}

\subsection{Lifting}

For $\alpha=(\alpha_1,\alpha_2)\in\mathbb R^2$ define
\begin{equation}\label{eq:lifting-symbol}
 m_\alpha(\xi,\eta)=|(\xi,\eta)|^{\alpha_1+\alpha_2}|\eta|^{\alpha_2},
\end{equation}
and let $\mathcal D^\alpha$ be the Fourier multiplier
$\widehat{\mathcal D^\alpha f}=m_\alpha\widehat f$.

\begin{theorem}[Lifting]\label{thm:lifting}
For $s,\alpha\in\mathbb R^2$, $0<p<\infty$ and
$0<q\leq\infty$, $\mathcal D^\alpha$ is an isomorphism
\begin{equation}\label{eq:lifting}
 \mathcal D^\alpha:\dot F^{s,q}_{p,F}
                   \longrightarrow\dot F^{s+\alpha,q}_{p,F},
 \qquad
 \|\mathcal D^\alpha f\|_{\dot F^{s+\alpha,q}_{p,F}}
                           \asymp\|f\|_{\dot F^{s,q}_{p,F}},
\end{equation}
with inverse $\mathcal D^{-\alpha}$; it restricts to an isomorphism
of the corresponding little spaces.
\end{theorem}

\subsection{Embeddings and Hardy-space consequences}

\begin{proposition}[Mixed-norm comparisons]\label{prop:embeddings}
Let $s\in\mathbb R^2$, $0<p<\infty$ and
$0<q\leq\infty$. Then
\begin{equation}\label{eq:FB-embeddings}
 \dot B^{s,\min(p,q)}_{p,F}
       \hookrightarrow\dot F^{s,q}_{p,F}
       \hookrightarrow\dot B^{s,\max(p,q)}_{p,F}.
\end{equation}
In particular $\dot F^{s,p}_{p,F}=\dot B^{s,p}_{p,F}$ with equality
of the displayed quasi-norms. For
$0<q_0\leq q_1\leq\infty$,
\begin{equation}\label{eq:q-monotonicity}
 \dot F^{s,q_0}_{p,F}\hookrightarrow\dot F^{s,q_1}_{p,F},
 \qquad
 \dot B^{s,q_0}_{p,F}\hookrightarrow\dot B^{s,q_1}_{p,F}.
\end{equation}
\end{proposition}

\begin{corollary}[Hardy-space consequences]\label{cor:Hardy-Besov}
For $0<p<\infty$,
\begin{equation}\label{eq:Hardy-Besov}
 \dot B^{0,\min(p,2)}_{p,F}
       \hookrightarrow H^p_F
       \hookrightarrow\dot B^{0,\max(p,2)}_{p,F}.
\end{equation}
If $1<p<\infty$ the middle space is $L^p$. For
$0<q_0\leq2\leq q_1\leq\infty$,
\begin{equation}\label{eq:Hardy-F}
 \dot F^{0,q_0}_{p,F}\hookrightarrow H^p_F
                         \hookrightarrow\dot F^{0,q_1}_{p,F}.
\end{equation}
\end{corollary}

The paper is organized as follows. Section~\ref{sec:discrete} proves
Theorem~\ref{thm:retraction}. Section~\ref{sec:tails} treats the
scale tails and proves Theorem~\ref{thm:little}.
Section~\ref{sec:carleson-duality} proves Theorems
\ref{thm:carleson-dual}--\ref{thm:Banach-endpoint} and
Corollary~\ref{cor:small-q-Besov}. Section~\ref{sec:full-dual}
proves Propositions \ref{prop:full-supremum} and
\ref{prop:modulations}. Section~\ref{sec:interpolation} proves
Theorems \ref{thm:complex-interpolation} and
\ref{thm:real-interpolation}. Section~\ref{sec:lifting} proves
Theorem~\ref{thm:lifting}. Section~\ref{sec:embeddings} proves
Proposition~\ref{prop:embeddings} and Corollary~\ref{cor:Hardy-Besov}.

\section{Discrete analysis and synthesis: proof of Theorem~\ref{thm:retraction}}
\label{sec:discrete}

Set $D=n+m$. We write $C$ for a positive constant, allowed to
change from one occurrence to the next and to depend on the
indices of the space and on finitely many seminorms of the
functions involved. Our Fourier convention is
$\widehat f(\xi,\eta)=\int e^{-i(x\cdot\xi+y\cdot\eta)}f(x,y)\,dx\,dy$.

The reproducing formulas, the flag grid and the pointwise sampling
estimates in this section are shared with the Besov framework of
\cite{LiBesov}. We give precise references for the common
identities and include the proofs of the sampling estimates used
below. The vector-valued maximal estimate of
Lemma~\ref{lem:maximal} and the $L^p(\ell^q)$ argument in
Lemma~\ref{lem:matrix} supply the bounds required for the
Triebel--Lizorkin scale.

The continuous Calder\'on reproducing formula
\begin{equation}\label{eq:continuousCRF}
 f=\sum_{j,k}\psi_{j,k}*\psi_{j,k}*f
\end{equation}
holds on $\mathcal S_F$ and on $\mathcal S_F'$: it is the continuous
formula of \cite{HL08,DLM10}, obtained by applying the two
partitions of unity \eqref{eq:partition1}--\eqref{eq:partition2}; the
mothers are real and even, flatness on $\eta=0$ gives the identity on
tests, and transposition gives it on distributions.

Sampling the continuous formula on the grids $\mathcal R_{j,k}^N$
gives the discrete Calder\'on reproducing formula, the discrete
version of \cite[Theorem~2.1]{HLW19}:
\begin{equation}\label{eq:exactCRF}
 f=\sum_{j,k}\sum_{R\in\mathcal R_{j,k}^N}
       |R|\,\psi_{j,k}(\,\cdot-z_R)\,(\psi_{j,k}*f)(z_R).
\end{equation}
At a fixed scale the identity is the Poisson summation formula on the
sampling lattice
$\Lambda_J=2^{-N}a_j\mathbb Z^n\times2^{-N}b_{j,k}\mathbb Z^m$;
the Fourier support box of $\psi_J*f$ is bounded by fixed multiples of
$a_j^{-1}$ and $b_{j,k}^{-1}$, so for a sufficiently large fixed $N$
no nonzero periodised copy meets it, and sampling reproduces the
$\psi_J$ piece. The cutoff and lattice estimates are identical in
notation to those of the companion paper
\cite[Section~2, (2.3)--(2.7)]{LiBesov}; we record the
convergence statement --- the series in \eqref{eq:exactCRF} converges
in $L^2$, in $\mathcal S_F$, and weakly in $\mathcal S_F'$ --- and use
the coefficient decay below.

For $\phi\in\mathcal S_F$, annular localisation and flatness give, for
arbitrary $L,M>0$,
\begin{equation}\label{eq:testcoeff}
 |(\psi_{j,k}*\phi)(z)|
 \leq C_{L,M,\phi}\,2^{-L(|j|+|k|)}(1+|z|)^{-M};
\end{equation}
this follows by differentiating $\widehat\psi_{j,k}\widehat\phi$, with
Taylor's formula at the low $\eta$ frequencies and the Schwartz bounds
at the high ones, and it is the estimate used throughout to justify
convergence of the synthesis series.

The same sampling gives the Plancherel identity: applying Poisson
summation to $|\psi_J*f|^2$ (whose Fourier support lies in the
difference of two support boxes, so the same $N$ applies) and then
Plancherel's theorem,
\begin{equation}\label{eq:L2isometry}
 \sum_{j,k}\sum_{R\in\mathcal R_{j,k}^N}|R|\,|(\psi_{j,k}*f)(z_R)|^2
 =\sum_{j,k}\|\psi_{j,k}*f\|_2^2=\|f\|_2^2,
\end{equation}
first for Schwartz functions and then on $L^2$ by bandlimited
approximation and monotone convergence.

The Plancherel--P\'olya comparison used below is the following:
for two admissible pairs $\psi,\varphi$ and sufficiently fine
fixed grids,
\begin{equation}\label{eq:PP}
 \bigl\|\{2^{-js_1}2^{-(j+k)s_2} U_{j,k}^\psi\}\bigr\|_{L^p(\ell^q)}
 \asymp
 \bigl\|\{2^{-js_1}2^{-(j+k)s_2} V_{j,k}^\varphi\}\bigr\|_{L^p(\ell^q)},
\end{equation}
where
\begin{equation}\label{eq:UV}
 U_{j,k}^{\psi}=\sum_R\sup_R|\psi_{j,k}*f|\,\mathbf1_R,\qquad
 V_{j,k}^{\varphi}=\sum_R\inf_R|\varphi_{j,k}*f|\,\mathbf1_R.
\end{equation}
We give a proof after the retraction estimates, using the sampling
estimates. The almost-orthogonality estimate used
in that proof and throughout the paper is
\begin{equation}\label{eq:almost}
 |\psi_{j,k}*\varphi_{j',k'}(z)|
 \leq C_{L,M}2^{-L(|j-j'|+|k-k'|)}P_{A,B,M}(z),
\end{equation}
where $A=\max(a_j,a_{j'})$, $B=\max(b_{j,k},b_{j',k'})$ and
\begin{equation}\label{eq:P}
 P_{A,B,M}(x,y)=
 \frac{A^M}{(A+|x|)^{n+M}}
 \frac{B^M}{(B+|y|)^{m+M}}.
\end{equation}
Indeed, the Fourier supports are disjoint unless both $|j-j'|$
and $|k-k'|$ are bounded by a fixed constant. For the remaining
pairs, integration by parts in the normalised frequency variables
gives arbitrary decay in $x/A$ and $y/B$, with constants uniform
over the nonzero kernels. This proves \eqref{eq:almost} for any
$L,M>0$.

We now prepare the maximal-function estimates that enter the
synthesis bound in the $L^p(\ell^q)$ ordering.

\begin{lemma}\label{lem:maximal}
For $1<u<\infty$ and $1<v\leq\infty$,
\begin{equation}\label{eq:vector-maximal}
 \left\|\left(\sum_\nu|M_Fh_\nu|^v\right)^{1/v}\right\|_{L^u}
 \leq C\left\|\left(\sum_\nu|h_\nu|^v\right)^{1/v}\right\|_{L^u}.
\end{equation}
The corresponding estimate with suprema holds for $v=\infty$.
\end{lemma}
\begin{proof}
Write $M_{n+m}$ for the ordinary maximal operator on
$\mathbb R^{n+m}$ and $M_m$ for the maximal operator in the $y$
variable. One has
\begin{equation}\label{eq:flag-domination}
                         M_Fh\leq C M_m(M_{n+m}h).
\end{equation}
Indeed, let $R=I\times J$ be a flag rectangle containing $(x,y)$,
with side lengths $a\leq b$. Put $J^*=J+[-a,a]^m$ and
$Q(v,a)=v+[-a,a]^m$. A cube in $\mathbb R^{n+m}$ of side comparable
to $a$, containing $(x,v)$, contains $I\times Q(v,a)$, so
\[
 M_{n+m}h(x,v)\geq c a^{-(n+m)}
                     \int_I\int_{Q(v,a)}|h(u,w)|\,dw\,du.
\]
For every $w\in J$, the set of $v\in J^*$ with $w\in Q(v,a)$ has
measure $(2a)^m$. Integrating in $v$, applying Fubini, and using
$|J^*|\asymp b^m$ and $y\in J^*$ therefore gives
\[
 M_m(M_{n+m}h)(x,y)
 \geq\frac1{|J^*|}\int_{J^*}M_{n+m}h(x,v)\,dv
 \geq\frac c{|R|}\int_R|h|.
\]
Taking the supremum over flag rectangles proves
\eqref{eq:flag-domination}. Two applications of the vector-valued
maximal inequality \cite{FS71}, with Fubini for the partial operator,
give, when $1<v<\infty$,
\begin{align}\label{eq:maximal-composition-expanded}
 \|\{M_Fh_\nu\}_\nu\|_{L^u(\ell^v)}
 &\leq C\|\{M_mM_{n+m}h_\nu\}_\nu\|_{L^u(\ell^v)}\notag\\
 &\leq C\|\{M_{n+m}h_\nu\}_\nu\|_{L^u(\ell^v)}
 \leq C\|\{h_\nu\}_\nu\|_{L^u(\ell^v)}.
\end{align}
At the supremum endpoint,
\begin{equation}\label{eq:maximal-supremum-expanded}
 \left\|\sup_\nu M_Fh_\nu\right\|_{L^u}
 \leq\left\|M_F\!\left(\sup_\nu|h_\nu|\right)\right\|_{L^u}
 \leq C\left\|\sup_\nu|h_\nu|\right\|_{L^u}.
\end{equation}
\end{proof}

\begin{lemma}\label{lem:grid}
Let $0<r\leq1$ and $M>D/r$. Choose arbitrary points
$z_{R'}\in R'\in\mathcal R_{j',k'}^N$ and nonnegative numbers
$c_{R'}$, and set $H=\sum_{R'}c_{R'}\mathbf1_{R'}$. If
$z,\bar z\in R\in\mathcal R_{j,k}^N$, then
\begin{align}\label{eq:grid-bound}
 \sum_{R'\in\mathcal R_{j',k'}^N} |R'|
          P_{A,B,M}(\bar z-z_{R'})c_{R'}
 \leq C\delta^{-D(1/r-1)}
       \left(\frac A{a'}\right)^{n(1/r-1)}
       \left(\frac B{b'}\right)^{m(1/r-1)}
                     [M_F(H^r)(z)]^{1/r},
\end{align}
where $\delta=2^{-N}$, $a'=a_{j'}$ and $b'=b_{j',k'}$. The two
scale-ratio factors are at most
$2^{D(1/r-1)(|j-j'|+|k-k'|)}$.
\end{lemma}
\begin{proof}
We recall the argument of \cite[Lemma~2.2]{LiBesov}.
It is enough to consider finitely many nonzero $c_{R'}$. Split
the points into product shells of radii $2^h A$ and $2^l B$,
$h,l\geq0$, with the zeroth shells including distances at most
$2A$, $2B$. On the $(h,l)$ shell the kernel is at most
$C A^{-n}B^{-m}2^{-h(n+M)}2^{-l(m+M)}$, and all source rectangles
in this shell, together with $z$, lie in a flag rectangle
$E_{h,l}$ with side lengths $C2^{\max(h,l)}A$ and
$C2^{\max(h,l)}B$. Since $A\leq B$, this common dilation preserves
the flag condition, and $|E_{h,l}|\leq C\,2^{D\max(h,l)}A^nB^m$.
At the source scale every rectangle has
volume $V=\delta^D(a')^n(b')^m$. Since $r\leq1$,
\begin{equation}\label{eq:shell-step}
 \sum_{R'\text{ in the shell}}Vc_{R'}
 \leq V^{1-1/r}
                  \left(\int_{E_{h,l}}H^r\right)^{1/r}.
\end{equation}
Combining these estimates, the shell contributes at most
\begin{align*}
 C\delta^{-D(1/r-1)}
 \left(\frac A{a'}\right)^{n(1/r-1)}
 \left(\frac B{b'}\right)^{m(1/r-1)}
 2^{-h(n+M)}2^{-l(m+M)}
 2^{D\max(h,l)/r}[M_F(H^r)(z)]^{1/r}.
\end{align*}
The required summability is explicit:
\begin{align}\label{eq:geometric-shell-sum}
 &\sum_{h,l\geq0}2^{-h(n+M)}2^{-l(m+M)}2^{D\max(h,l)/r}\notag\\
 &\quad\leq\left(\sum_{h\geq0}2^{-h(n+M-D/r)}\right)
             \left(\sum_{l\geq0}2^{-l(m+M-D/r)}\right)\notag\\
 &\quad=\frac1{(1-2^{-(n+M-D/r)})(1-2^{-(m+M-D/r)})}<\infty.
\end{align} The last assertion follows from
$\log_2(A/a')=(j'-j)_+$ and
$\log_2(B/b')=(\min(j',k')-\min(j,k))_+$, each bounded by
$|j-j'|+|k-k'|$.
The argument also permits different fixed source and target
grids: the factor $\delta^{-D(1/r-1)}$ then uses the source
mesh, and the constant may depend on both fixed meshes.
\end{proof}

\begin{lemma}\label{lem:matrix}
Suppose nonnegative functions $u_{j,k}$ satisfy, for
$z\in R\in\mathcal R_{j,k}^N$,
\begin{align}\label{eq:matrix-hypothesis}
 u_{j,k}(z)\leq C\sup_{\bar z\in R}\sum_{j',k'}
   2^{-L(|j-j'|+|k-k'|)}
   \sum_{R'\in\mathcal R_{j',k'}^N}|R'|
          P_{A,B,M}(\bar z-z_{R'})c_{j',k',R'}.
\end{align}
Set $H_{j,k}=\sum_R c_{j,k,R}\mathbf1_R$. Let $0<p<\infty$,
$0<q\leq\infty$, and choose
$0<r<\min(1,p,q)$ when $q<\infty$ and
$0<r<\min(1,p)$ when $q=\infty$. If
$M>D/r$ and $L>|s_1+s_2|+|s_2|+D(1/r-1)$, then
\begin{equation}\label{eq:matrix-F}
 \|\{2^{-js_1}2^{-(j+k)s_2} u_{j,k}\}\|_{L^p(\ell^q)}
 \leq C\|\{2^{-js_1}2^{-(j+k)s_2} H_{j,k}\}\|_{L^p(\ell^q)}.
\end{equation}
\end{lemma}
\begin{proof}
The weight ratio satisfies
\begin{equation}\label{eq:smoothness-transfer}
 2^{-js_1}2^{-(j+k)s_2}2^{j's_1}2^{(j'+k')s_2}
 \leq 2^{|s_1+s_2||j-j'|+|s_2||k-k'|}.
\end{equation}
Lemma~\ref{lem:grid} therefore gives
\begin{equation}\label{eq:pointwise-matrix}
 2^{-js_1}2^{-(j+k)s_2} u_{j,k}
 \leq C\sum_{j',k'}2^{-\varepsilon(|j-j'|+|k-k'|)}
       [M_F((2^{-j's_1}2^{-(j'+k')s_2}H_{j',k'})^r)]^{1/r},
\end{equation}
with $\varepsilon=L-|s_1+s_2|-|s_2|-D(1/r-1)>0$. Write
$a_{u,v}=2^{-\varepsilon(|u|+|v|)}$. For a nonnegative sequence $b$,
\begin{align}\label{eq:sequence-convolution-expanded}
 \|a*b\|_{\ell^q}
 &\leq\|a\|_{\ell^1}\|b\|_{\ell^q},&&1\leq q\leq\infty,\notag\\
 \|a*b\|_{\ell^q}^q
 &\leq\sum_{J,L}a_{J-L}^qb_L^q
 =\|a\|_{\ell^q}^q\|b\|_{\ell^q}^q,&&0<q<1.
\end{align}
Both constants are finite because
\[
 \sum_{u,v\in\mathbb Z}2^{-\varepsilon\rho(|u|+|v|)}
 =\left(\frac{1+2^{-\varepsilon\rho}}
              {1-2^{-\varepsilon\rho}}\right)^2<\infty
 \qquad(\rho>0).
\]
Set $B_J=2^{-js_1}2^{-(j+k)s_2}H_J$. The remaining mixed norm obeys
\begin{align}\label{eq:matrix-maximal-expanded}
 \|\{[M_F(B_J^r)]^{1/r}\}_J\|_{L^p(\ell^q)}^r
 &=\|\{M_F(B_J^r)\}_J\|_{L^{p/r}(\ell^{q/r})}\notag\\
 &\leq C\|\{B_J^r\}_J\|_{L^{p/r}(\ell^{q/r})}
 =C\|\{B_J\}_J\|_{L^p(\ell^q)}^r.
\end{align}
Here $p/r>1$ and $q/r>1$, as required by
Lemma~\ref{lem:maximal}. For $q=\infty$,
\begin{align}\label{eq:matrix-maximal-sup-expanded}
 \left\|\sup_J[M_F(B_J^r)]^{1/r}\right\|_{L^p}^r
 &\leq\left\|M_F\!\left(\sup_J B_J^r\right)\right\|_{L^{p/r}}\notag\\
 &\leq C\left\|\sup_J B_J^r\right\|_{L^{p/r}}.
\end{align}
Together with \eqref{eq:pointwise-matrix}, these estimates prove
\eqref{eq:matrix-F}.
\end{proof}

\begin{lemma}\label{lem:peetre}
Let $0<a\leq b$, $0<r\leq1$ and $\lambda r>D$. If the Fourier support of a
tempered distribution $h$ is contained in
$\{(\xi,\eta):|\xi|\leq C_0/a,\ |\eta|\leq C_0/b\}$, then
\begin{equation}\label{eq:regularized-peetre}
 h^*_{a,b,\lambda}(z)
       \leq C[M_F(|h|^r)(z)]^{1/r},
\end{equation}
where
$\rho_{a,b}(x,y)=\max(|x|/a,|y|/b)$ and
$h^*_{a,b,\lambda}(z)=\sup_u|h(u)|(1+\rho_{a,b}(z-u))^{-\lambda}$.
For $u,z$ in the same rectangle of side lengths $\delta a,\delta b$,
$0<\delta\leq1$,
\begin{equation}\label{eq:normalised-derivative}
 a|\nabla_x h(u)|+b|\nabla_y h(u)|
       \leq C[M_F(|h|^r)(z)]^{1/r}.
\end{equation}
\end{lemma}
\begin{proof}
This estimate is also proved in \cite[Lemma~2.4]{LiBesov}; we give
the argument here.
We first suppose that $h\in\mathcal S$, so that
$P=h^*_{a,b,\lambda}(z)<\infty$. Choose
$\gamma\in\mathcal S$ with $\widehat\gamma\in C_c^\infty$ equal to one
on a neighbourhood of the normalised support box, and put
$\gamma_{a,b}(x,y)=a^{-n}b^{-m}\gamma(x/a,y/b)$. Then
$h=\gamma_{a,b}*h$. Write $\rho=\rho_{a,b}$. The inequalities
\[
 |h(v)|^{1-r}\leq P^{1-r}(1+\rho(z-v))^{\lambda(1-r)},
 \qquad
 1+\rho(z-v)\leq(1+\rho(z-u))(1+\rho(u-v))
\]
and the Schwartz decay of $\gamma$ give, for any $T>D$,
\begin{equation}\label{eq:peetre-step}
 \frac{|h(u)|}{(1+\rho(z-u))^\lambda}
 \leq C P^{1-r}(1+\rho(z-u))^{-\lambda r}
 \int a^{-n}b^{-m}(1+\rho(u-v))^{-T}|h(v)|^r\,dv.
\end{equation}
Here the zero function is treated separately. Put $d=\rho(z-u)$.
For each $l\geq0$, the set $\{v:\rho(u-v)\leq2^{l+1}\}$ and the
point $z$ are contained in a flag rectangle $E_l$ with
$|E_l|\leq C a^nb^m(2^l+d)^D$. A dyadic shell decomposition thus yields
\begin{align}\label{eq:peetre-shells-expanded}
 &\int a^{-n}b^{-m}(1+\rho(u-v))^{-T}|h(v)|^r\,dv\notag\\
 &\qquad\leq C\sum_{l\geq0}2^{-lT}(2^l+d)^D
                      M_F(|h|^r)(z)\notag\\
 &\qquad\leq C(1+d)^D\sum_{l\geq0}2^{-l(T-D)}
                      M_F(|h|^r)(z)
 \leq C(1+d)^D M_F(|h|^r)(z).
\end{align}
Combining this with \eqref{eq:peetre-step}, and using
$D-\lambda r<0$, gives
\[
 P\leq C P^{1-r}M_F(|h|^r)(z),
 \qquad P^r\leq C M_F(|h|^r)(z),
\]
which proves \eqref{eq:regularized-peetre}.
The normalised derivatives $a\nabla_x\gamma_{a,b}$ and
$b\nabla_y\gamma_{a,b}$ satisfy the same Schwartz kernel bounds.
Differentiating the reproducing identity and using
$|h(v)|\leq P(1+\rho(z-v))^\lambda$ therefore gives
\[
 a|\nabla_xh(u)|+b|\nabla_yh(u)|
 \leq C P(1+\rho(z-u))^\lambda.
\]
For $u,z$ in one of the stated rectangles, $\rho(z-u)\leq C\delta\leq C$;
hence \eqref{eq:normalised-derivative} follows from
\eqref{eq:regularized-peetre}.

For a general band-limited tempered distribution $h$, its smooth
representative and all its derivatives have polynomial growth.
Choose $\vartheta\in\mathcal S$ whose Fourier transform is a
nonnegative function in $C_c^\infty$, normalised so that
$\vartheta(0)=1$. Then $|\vartheta|\leq1$. Set
\[
 h_\epsilon(x,y)=h(x,y)\vartheta(\epsilon x/a,\epsilon y/b),
 \qquad 0<\epsilon\leq1.
\]
Each $h_\epsilon$ is Schwartz, and the Fourier supports of all these
functions lie in one box of the stated form with a fixed larger
constant in place of $C_0$. The preceding estimates consequently
hold with constants independent of $\epsilon$. Since
$|h_\epsilon|\leq|h|$ and $h_\epsilon$, together with its first
derivatives, converges pointwise to $h$, letting
$\epsilon\downarrow0$ proves both assertions without an a priori
finiteness assumption on the Peetre supremum.
\end{proof}

\begin{lemma}\label{lem:Linfty-sampling}
For a sufficiently fine fixed grid,
\begin{equation}\label{eq:Linfty-sampling}
 \|\psi_J*g\|_{L^\infty}
       \asymp\sup_{R\in\mathcal R_J^N}|(\psi_J*g)(z_R)|,
\end{equation}
uniformly in $J$ and $g\in\mathcal S_F'$.
\end{lemma}
\begin{proof}
This estimate is also proved in \cite[Lemma~2.5]{LiBesov}; we give
the argument here.
The sampled supremum is bounded by the continuous one. Write
$h=\psi_J*g$, $a=a_j$ and $b=b_{j,k}$, and first assume that $h$ is
bounded. By \eqref{eq:normalised-derivative} with $r=1$,
\[
 a\|\nabla_xh\|_\infty+b\|\nabla_yh\|_\infty
 \leq C\|h\|_\infty.
\]
For $u\in R$, the mean value theorem along the segment from $z_R$
to $u$ gives
\[
 |h(u)|\leq |h(z_R)|+C2^{-N}\|h\|_\infty.
\]
Taking the supremum and choosing $N$ such that $C2^{-N}\leq1/2$
yields
\[
 \|h\|_\infty\leq2\sup_{R\in\mathcal R_J^N}|h(z_R)|.
\]
For an arbitrary $h=\psi_J*g$, apply the bounded case to the
regularisations $h_\epsilon$ from the proof of
Lemma~\ref{lem:peetre}. Their Fourier supports lie in one common
enlarged box, so the same sufficiently large $N$ works for every
$\epsilon$. Moreover,
\[
 \sup_R|h_\epsilon(z_R)|\leq\sup_R|h(z_R)|.
\]
Letting $\epsilon\downarrow0$ pointwise and then taking the
supremum proves the reverse inequality in full generality.
\end{proof}

We are now in a position to prove Theorem~\ref{thm:retraction}.
For the analysis bound choose $r$ as in Lemma~\ref{lem:matrix}. On every rectangle,
Lemma~\ref{lem:peetre} gives
\[
 2^{-js_1}2^{-(j+k)s_2}\sup_R|\psi_{j,k}*f|
 \leq C\bigg(M_F(|2^{-js_1}2^{-(j+k)s_2}(\psi_{j,k}*f)|^r)(z)\bigg)^{1/r},
 \qquad z\in R.
\]
Taking the mixed norm and using
\eqref{eq:matrix-maximal-expanded} or
\eqref{eq:matrix-maximal-sup-expanded} gives
\begin{align}\label{eq:analysis-chain-expanded}
 \|S_\psi f\|_{\dot f^{s,q}_{p,F}}
 &\leq\|\{2^{-js_1}2^{-(j+k)s_2}U_J^\psi\}_J\|_{L^p(\ell^q)}\notag\\
 &\leq C\|\{[M_F(|2^{-js_1}2^{-(j+k)s_2}\psi_J*f|^r)]^{1/r}\}_J
                    \|_{L^p(\ell^q)}\notag\\
 &\leq C\|f\|_{\dot F^{s,q}_{p,F}}.
\end{align}
We note that this analysis bound does not use \eqref{eq:PP}.
For the synthesis estimate, a coefficient sequence satisfies
\begin{equation}\label{eq:coefficient-growth}
 |t_{j,k,R}|
 \leq (2^{-js_1}2^{-(j+k)s_2})^{-1}|R|^{-1/p}
                         \|t\|_{\dot f^{s,q}_{p,F}},
\end{equation}
and $|R|^{-1/p}\leq C2^{C_1(|j|+|k|)}$. For $M>D$ the lattice
counting bound
\begin{equation}\label{eq:lattice-counting}
 \sum_{R\in\mathcal R_{j,k}^N}|R|(1+|z_R|)^{-M}
 \leq C_M(1+2^{-N}a_j)^n(1+2^{-N}b_{j,k})^m
\end{equation}
holds. In view of \eqref{eq:testcoeff}, the series
\begin{equation}\label{eq:synthesis-abs}
 \sum_{j,k,R}|R|\,|t_{j,k,R}|\,
                  |\langle\psi_{j,k}(\,\cdot-z_R),h\rangle|
\end{equation}
converges for every $h\in\mathcal S_F$. Indeed, the smoothness and
cell-volume factors grow at most exponentially in $|j|+|k|$, so for
a fixed $C_2$ and then $L>C_2+1$,
\begin{align}\label{eq:distribution-convergence-expanded}
 &\sum_{j,k,R}|R|\,|t_{j,k,R}|
           |\langle\psi_{j,k}(\,\cdot-z_R),h\rangle|\notag\\
 &\qquad\leq C\|t\|_{\dot f^{s,q}_{p,F}}\mathfrak p(h)
       \sum_{j,k}2^{-(L-C_2)(|j|+|k|)}
 \leq C\mathfrak p(h)\|t\|_{\dot f^{s,q}_{p,F}}.
\end{align}
Thus $T_\psi t$ is continuous on $\mathcal S_F$.
For finite $t$, convolution of synthesis with $\varphi_J$ gives
\begin{align}\label{eq:synthesis-pointwise-expanded}
 |\varphi_J*T_\psi t(z)|
 &\leq\sum_L\sum_{R\in\mathcal R_L^N}|R|\,|t_{L,R}|
                         |(\varphi_J*\psi_L)(z-z_R)|\notag\\
 &\leq C\sum_L2^{-L_0(|j-l_1|+|k-l_2|)}
      \sum_{R\in\mathcal R_L^N}|R|P_{A,B,M}(z-z_R)|t_{L,R}|.
\end{align}
Choose $L_0>|s_1+s_2|+|s_2|+D(1/r-1)$ and $M>D/r$.
Lemma~\ref{lem:matrix} then gives
\begin{equation}\label{eq:synthesis-bound}
 \|T_\psi t\|_{\dot F^{s,q}_{p,F}}
                  \leq C\|t\|_{\dot f^{s,q}_{p,F}}.
\end{equation}
Finite truncations are uniformly bounded and converge in
$\mathcal S_F'$ by the absolute convergence; Fatou's lemma gives
the same estimate for arbitrary sequences.

If $p,q<\infty$, let $\Gamma_a$ be an increasing finite exhaustion
of the labels. The coefficient tail satisfies
\begin{align}\label{eq:coefficient-density-expanded}
 G_a(z)^q
 &:=\sum_{(J,R)\notin\Gamma_a}
   [2^{-js_1}2^{-(j+k)s_2}|t_{J,R}|]^q\mathbf1_R(z)
 \downarrow0\quad\text{a.e.},\notag\\
 G_a(z)^p&\leq
 \left(\sum_{J,R}[2^{-js_1}2^{-(j+k)s_2}|t_{J,R}|]^q
                         \mathbf1_R(z)\right)^{p/q}\in L^1,\notag\\
 \|t-\mathbf1_{\Gamma_a}t\|_{\dot f^{s,q}_{p,F}}^p
 &=\int G_a^p\longrightarrow0.
\end{align}
This proves density of finite sequences.
The synthesis bound then gives norm convergence of the series, and
the weak convergence of every synthesis series follows from
\eqref{eq:synthesis-abs}. The identity $T_\psi S_\psi=I$ is
\eqref{eq:exactCRF}; taking $t=S_\psi f$ and using the synthesis
estimate proves \eqref{eq:discrete-norm}.

It remains to record completeness and the basic properties. The
coefficient space $\dot f^{s,q}_{p,F}$ is complete: after
multiplication by $2^{-js_1}2^{-(j+k)s_2}$ it is a closed subspace of
$L^p(\ell^q)$, and $S_\psi T_\psi$ is a bounded projection whose
range is isomorphic to $\dot F^{s,q}_{p,F}$ by
\eqref{eq:discrete-norm}. For $p,q\geq1$ the continuous expression
is a norm, by the two Minkowski inequalities and definiteness from
the reproducing formula. The estimate
\eqref{eq:coefficient-growth} with $t=S_\psi f$ proves the
continuous embedding into $\mathcal S_F'$. Every test function
belongs to the space by \eqref{eq:testcoeff} and the resulting
spatial decay, and finite synthesis sums are test functions; the
density for $q<\infty$ and the weak approximation for $q=\infty$
follow from the density of finite sequences and the bounded
synthesis. This proves all the assertions of
Theorem~\ref{thm:retraction}.

For completeness, we now prove \eqref{eq:PP}. Choose
$0<r<\min(1,p,q)$, with $q$ omitted when $q=\infty$. First use
the same analysing pair on both sides and write $h_{j,k}=\psi_{j,k}*f$,
$U_{j,k}=U_{j,k}^{\psi}$ and $V_{j,k}=V_{j,k}^{\psi}$. The mean
value theorem and \eqref{eq:normalised-derivative} imply, on each
rectangle,
\[
 U_{j,k}\leq V_{j,k}
       +C2^{-N}[M_F(|h_{j,k}|^r)]^{1/r}.
\]
Put $\tau=\min(1,p,q)$, with $q$ omitted when $q=\infty$.
The $\tau$-triangle inequality and Lemma~\ref{lem:maximal} give
\[
\begin{aligned}
 \|\{2^{-js_1}2^{-(j+k)s_2}U_{j,k}\}\|_{L^p(\ell^q)}^{\tau}
 &\leq \|\{2^{-js_1}2^{-(j+k)s_2}V_{j,k}\}\|_{L^p(\ell^q)}^{\tau}\\
 &\quad +(C2^{-N})^{\tau}
      \|\{2^{-js_1}2^{-(j+k)s_2}U_{j,k}\}\|_{L^p(\ell^q)}^{\tau}.
\end{aligned}
\]
Choose the common grid resolution so that
$\kappa=(C2^{-N})^\tau<1/2$. Absorption gives
\begin{equation}\label{eq:PP-absorption-expanded}
 \|\{2^{-js_1}2^{-(j+k)s_2}U_{j,k}\}\|_{L^p(\ell^q)}
 \leq(1-\kappa)^{-1/\tau}
       \|\{2^{-js_1}2^{-(j+k)s_2}V_{j,k}\}\|_{L^p(\ell^q)}.
\end{equation}
To justify the argument without an a priori finiteness assumption,
first restrict to finitely many scales and replace each $h_{j,k}$
by its bandlimited regularisation from Lemma~\ref{lem:peetre}.
The regularised rectangle infima are bounded by $V_{j,k}$ because
$|\vartheta|\leq1$. Absorption is uniform for the common enlarged
Fourier support box. Pointwise convergence on each rectangle and
Fatou's lemma, followed by exhaustion of the scales, give the
assertion for every distribution. Finally, synthesis with one
admissible pair and analysis with the other, as proved above,
compare their continuous quasi-norms. Combining that comparison
with the same-pair supremum--infimum estimate proves
\eqref{eq:PP} and the independence of the analysing pair.

\section{Scale tails and the little space: proof of Theorem~\ref{thm:little}}
\label{sec:tails}

Fix an admissible pair and a sufficiently fine grid. We first record
an estimate that permits restriction to an arbitrary set of scales.

\begin{lemma}[Restricted sampling]\label{lem:restricted-sampling}
For every set $E\subset\mathbb Z^2$ and $f\in\mathcal S_F'$,
\begin{equation}\label{eq:restricted-sampling}
 \bigl\|\{2^{-js_1}2^{-(j+k)s_2} U_J\}_{J\in E}\bigr\|_{L^p(\ell^q)}
 \leq C\bigl\|\{2^{-js_1}2^{-(j+k)s_2}|\psi_J*f|\}_{J\in E}\bigr\|_{L^p(\ell^q)},
\end{equation}
where $U_J=\sum_R\sup_R|\psi_J*f|\,\mathbf1_R$. The constant is
independent of $E$ and $f$; the same estimate holds with the rectangle
suprema replaced by the sampled values $|t_{J,R}|$.
\end{lemma}
\begin{proof}
Choose $r$ as in Lemma~\ref{lem:matrix} and $\lambda$ with
$\lambda r>D$. On each grid rectangle the regularised Peetre estimate
\eqref{eq:regularized-peetre} gives
$2^{-js_1}2^{-(j+k)s_2} U_J(z)\leq C[M_F(|2^{-js_1}2^{-(j+k)s_2}(\psi_J*f)|^r)(z)]^{1/r}$. For finite
$q$ apply Lemma~\ref{lem:maximal} with exponents $p/r,q/r$ and take
the $r$-th root; for $q=\infty$ use the pointwise supremum and the
scalar maximal estimate. The sampled values are bounded by the rectangle
suprema.
\end{proof}

For $K\geq0$ and $L\geq1$ define the finite synthesis
\begin{equation}\label{eq:finite-sum}
 f_{K,L}
 =\sum_{|J|_1\leq K}
   \sum_{\substack{R\in\mathcal R_{j,k}^N\\|z_R|\leq L}}
       |R|\,t_{J,R}\psi_{j,k}(\,\cdot-z_R),
\end{equation}
and the sampled spatial tail
\begin{equation}\label{eq:space-tail}
 \tau^F_{f,q}(K,L)
 =\left\|\left\{
      2^{-js_1}2^{-(j+k)s_2}\sum_{\substack{R\in\mathcal R_{j,k}^N\\|z_R|>L}}
                      |t_{J,R}|\mathbf1_R
               \right\}_{|J|_1\leq K}\right\|_{L^p(\ell^q)}.
\end{equation}
Every finite sum $f_{K,L}$ belongs to $\mathcal S_F$.

\begin{proposition}[Quantitative approximation]\label{prop:quantitative}
For $f\in\dot F^{s,q}_{p,F}$,
\begin{align}
 \|f-f_{K,L}\|_{\dot F^{s,q}_{p,F}}
 &\leq C\bigl(\omega^F_{f,q}(K)+\tau^F_{f,q}(K,L)\bigr),
                    \label{eq:approximation-error}\\
 \|f_{K,L}\|_{\dot F^{s,q}_{p,F}}
 &\leq C\|f\|_{\dot F^{s,q}_{p,F}}.\label{eq:approximation-uniform}
\end{align}
The constants are independent of $K,L,f$. For fixed $K$,
$\tau^F_{f,q}(K,L)\to0$ as $L\to\infty$. The sums $f_{K,L}$
converge weakly to $f$ in $\mathcal S_F'$; if $q<\infty$ they give
norm approximation by first choosing $K$ and then $L$.
\end{proposition}
\begin{proof}
Let $t=S_\psi f$ and let $E_{K,L}$ be the finite set of retained
labels. Split the omitted coefficients as
\begin{align}\label{eq:tail-decomposition-expanded}
 t^{\rm sc}_{J,R}&=t_{J,R}\mathbf1_{\{|J|_1>K\}},\notag\\
 t^{\rm sp}_{J,R}&=t_{J,R}
       \mathbf1_{\{|J|_1\leq K,\ |z_R|>L\}},\notag\\
 f-f_{K,L}&=T_\psi t^{\rm sc}+T_\psi t^{\rm sp},\qquad
 f_{K,L}=T_\psi(\mathbf1_{E_{K,L}}t).
\end{align}
Restricted sampling and the definition of the spatial tail give
\[
 \|t^{\rm sc}\|_{\dot f^{s,q}_{p,F}}\leq C\omega^F_{f,q}(K),
 \qquad
 \|t^{\rm sp}\|_{\dot f^{s,q}_{p,F}}=\tau^F_{f,q}(K,L).
\]
With $\tau=\min(1,p,q)$, omitting $q$ when $q=\infty$, the bounded
synthesis map gives
\begin{align}\label{eq:tail-error-expanded}
 \|f-f_{K,L}\|_{\dot F^{s,q}_{p,F}}^\tau
 &\leq C\bigl(\|t^{\rm sc}\|_{\dot f^{s,q}_{p,F}}^\tau
              +\|t^{\rm sp}\|_{\dot f^{s,q}_{p,F}}^\tau\bigr)\notag\\
 &\leq C\bigl(\omega^F_{f,q}(K)^\tau
                  +\tau^F_{f,q}(K,L)^\tau\bigr),\\
 \|f_{K,L}\|_{\dot F^{s,q}_{p,F}}
 &\leq C\|\mathbf1_{E_{K,L}}t\|_{\dot f^{s,q}_{p,F}}
 \leq C\|t\|_{\dot f^{s,q}_{p,F}}
 \leq C\|f\|_{\dot F^{s,q}_{p,F}}.\label{eq:finite-bound-expanded}
\end{align}
The inequality $(a^\tau+b^\tau)^{1/\tau}\leq2^{1/\tau-1}(a+b)$
proves \eqref{eq:approximation-error}.

For fixed $K$ only finitely many scales occur, and at almost every
point the spatial tail tends to zero as $L\to\infty$; it is
dominated by the full sampled expression, which lies in $L^p$ by
bounded analysis. Dominated convergence proves the claim. The
weak convergence follows from the absolute convergence of
synthesis, and for $q<\infty$ the pointwise scale tails decrease
to zero and are dominated by the full defining expression, so
dominated convergence gives $\omega^F_{f,q}(K)\to0$ and norm
approximation.
\end{proof}

We are now ready to prove Theorem~\ref{thm:little}. The tail
$\omega^F_{f,\infty}(K)$ decreases with $K$ and is bounded by the
defining norm, so $d_F(f)$ exists. Every $g\in\mathcal S_F$
satisfies $\omega^F_{g,\infty}(K)\to0$ by \eqref{eq:testcoeff}. Put
$\tau=\min(1,p)$. The pointwise supremum inequality and the
$\tau$-triangle inequality in $L^p$ give
\begin{equation}\label{eq:little-lower}
 [\omega^F_{f,\infty}(K)]^\tau
 \leq\|f-g\|_{\dot F^{s,\infty}_{p,F}}^\tau
                    +[\omega^F_{g,\infty}(K)]^\tau.
\end{equation}
Letting $K\to\infty$ and taking the infimum over $g$ proves the
lower bound in \eqref{eq:distance}. For the upper bound, \eqref{eq:approximation-error} gives
\begin{align}\label{eq:distance-upper-expanded}
 \operatorname{dist}(f,\mathcal S_F)
 &\leq\liminf_{L\to\infty}\|f-f_{K,L}\|_{\dot F^{s,\infty}_{p,F}}\notag\\
 &\leq C\omega^F_{f,\infty}(K)
 \longrightarrow C d_F(f)\qquad(K\to\infty).
\end{align} The
vanishing-tail condition is thus equivalent to membership in the
test-function closure, and when it holds the same choices give
approximation by finite syntheses; conversely every finite synthesis
is a test function, so the two closures coincide.

For the pointwise form, the decreasing functions
$H_K(z)=\sup_{|j|+|k|>K}2^{-js_1}2^{-(j+k)s_2}|\psi_{j,k}*f(z)|$ are dominated
by a function in $L^p$. If $H_\infty=\lim_KH_K$, dominated
convergence gives
\begin{equation}\label{eq:tail-limit-expanded}
 d_F(f)^p=\lim_{K\to\infty}\int H_K(z)^p\,dz
           =\int H_\infty(z)^p\,dz.
\end{equation} Hence $d_F(f)=0$ exactly when
$H_K(z)\to0$ almost everywhere, which is
\eqref{eq:pointwise-little}. Equivalent admissible-pair quasi-norms
have the same closure of $\mathcal S_F$, and the continuous
condition \eqref{eq:little} contains no grid, so the characterisation
is independent of the pair and the grid.

\begin{remark}[Approximation rates]\label{rem:rate}
The estimate \eqref{eq:approximation-error} records the actual
scale and spatial tails; membership in the little space alone
prescribes no rate. If $\omega^F_{f,q}(K)\leq C_0 2^{-\epsilon K}$
and, for some $\gamma>0$,
\begin{equation}\label{eq:rate-hypothesis}
 \left\|\left\{2^{-js_1}2^{-(j+k)s_2}\sum_R
       (1+|z_R|)^\gamma|t_{J,R}|\mathbf1_R\right\}_J
       \right\|_{L^p(\ell^q)}\leq C_1,
\end{equation}
then
\begin{equation}\label{eq:rate-conclusion}
 \|f-f_{K,L}\|_{\dot F^{s,q}_{p,F}}
 \leq C\bigl(C_0 2^{-\epsilon K}+C_1(1+L)^{-\gamma}\bigr).
\end{equation}
\end{remark}

\section{Carleson duality: proofs of Theorems~\ref{thm:carleson-dual}--\ref{thm:Banach-endpoint} and Corollary~\ref{cor:small-q-Besov}}
\label{sec:carleson-duality}

For sequences we use the pairing
\begin{equation}\label{eq:sequence-pairing}
 [t,d]=\sum_{j,k}\sum_R|R|\,t_{j,k,R}\overline{d_{j,k,R}}.
\end{equation}

\begin{lemma}[Carleson sequence duality]\label{lem:carleson-sequence}
Under the assumptions of Theorem~\ref{thm:carleson-dual},
\begin{equation}\label{eq:carleson-sequence}
 \bigl(\dot f^{s,q}_{p,F}\bigr)^*=\dot c_F^{-s,p,q'}.
\end{equation}
The pairing converges absolutely and the two norms are equivalent.
\end{lemma}
\begin{proof}
First let $t$ have finite support and put
$G(z)=\{\sum_{j,k,R}[2^{-js_1}2^{-(j+k)s_2}|t_{j,k,R}|]^q\mathbf1_R(z)\}^{1/q}$
and $\Omega_l=\{z:G(z)>2^l\}$. Assign a nonzero coefficient to
$\mathcal B_l$ when
$|R\cap\Omega_l|>\tfrac12|R|$ and
$|R\cap\Omega_{l+1}|\leq\tfrac12|R|$. Every nonzero coefficient has
exactly one such index. Let
$\Omega_l^*=\{z:M_F\mathbf1_{\Omega_l}(z)>1/4\}$. This is open,
contains every rectangle in $\mathcal B_l$, and the scalar maximal
inequality gives
\begin{equation}\label{eq:carleson-enlargement}
 |\Omega_l^*|\leq C|\Omega_l|.
\end{equation}
For $R\in\mathcal B_l$, the set $E_R=R\setminus\Omega_{l+1}$ has
$|E_R|\geq|R|/2$, so
\begin{equation}\label{eq:carleson-group}
 \sum_{\mathcal B_l}[2^{-js_1}2^{-(j+k)s_2}|t_{j,k,R}|]^q|R|
 \leq C\int_{\Omega_l^*\setminus\Omega_{l+1}}G(z)^q\,dz
 \leq C2^{lq}|\Omega_l^*|.
\end{equation}
H\"older's inequality for this group gives
\begin{align}
 \sum_{\mathcal B_l}|R|\,|t_{j,k,R}d_{j,k,R}|
 &\leq C2^l|\Omega_l^*|^{1/q}
   \left\{\sum_{R\subset\Omega_l^*}|R|
   [(2^{-js_1}2^{-(j+k)s_2})^{-1}|d_{j,k,R}|]^{q'}\right\}^{1/q'}\notag\\
 &\leq C\|d\|_{\dot c_F^{-s,p,q'}}2^l|\Omega_l^*|^{1/p}.
\end{align}
For $p\leq1$, the scale summation and the dyadic layer-cake
estimate are
\begin{align}\label{eq:carleson-layercake-expanded}
 \sum_l2^l|\Omega_l^*|^{1/p}
 &\leq\left(\sum_l2^{lp}|\Omega_l^*|\right)^{1/p}
 \leq C\left(\sum_l2^{lp}|\Omega_l|\right)^{1/p},\notag\\
 \sum_l2^{lp}|\Omega_l|
 &=\int\sum_{\{l:\,2^l<G(z)\}}2^{lp}\,dz
 \leq\frac{1}{1-2^{-p}}\int G(z)^p\,dz.
\end{align}
Combining these estimates with the preceding group bound gives
\begin{equation}\label{eq:carleson-pairing}
 \sum_{j,k,R}|R|\,|t_{j,k,R}d_{j,k,R}|
 \leq C\|d\|_{\dot c_F^{-s,p,q'}}
   \{\sum_l2^{lp}|\Omega_l|\}^{1/p}
 \leq C\|d\|_{\dot c_F^{-s,p,q'}}\|t\|_{\dot f^{s,q}_{p,F}}.
\end{equation}
Finite restrictions and monotone convergence give absolute
convergence for arbitrary $t$.

Conversely, let $\lambda$ be a bounded functional and define $d$ by
$\lambda(e_{j,k,R})=|R|\overline{d_{j,k,R}}$. Fix an open $\Omega$
of finite positive measure and a finite collection $\Gamma$ of
labelled rectangles contained in it. Since $p<q$, a sequence supported on
this collection satisfies
\begin{equation}\label{eq:carleson-test}
 \|t\|_{\dot f^{s,q}_{p,F}}
 \leq|\Omega|^{1/p-1/q}
   \{\sum[2^{-js_1}2^{-(j+k)s_2}|t_{j,k,R}|]^q|R|\}^{1/q}.
\end{equation}
The finite-dimensional duality step can be written explicitly.
Put
\[
 A_{J,R}=2^{-js_1}2^{-(j+k)s_2},\qquad
 B_{J,R}=A_{J,R}^{-1}|d_{J,R}|,\qquad
 B_\Gamma=\left(\sum_{\Gamma}|R|B_{J,R}^{q'}\right)^{1/q'}.
\]
If $B_\Gamma>0$, choose the finite sequence
\begin{equation}\label{eq:carleson-normer-expanded}
 t_{J,R}=A_{J,R}^{-1}
       \frac{B_{J,R}^{q'-1}}{B_\Gamma^{q'/q}}\,
       \frac{d_{J,R}}{|d_{J,R}|}\,\mathbf1_\Gamma(J,R),
\end{equation}
with value zero when $d_{J,R}=0$. Since $(q'-1)q=q'$,
\begin{align}\label{eq:carleson-normer-check-expanded}
 \sum_{\Gamma}|R|(A_{J,R}|t_{J,R}|)^q&=1,\notag\\
 \lambda(t)=\sum_{\Gamma}|R|t_{J,R}\overline{d_{J,R}}
 &=B_\Gamma^{\,q'-q'/q}=B_\Gamma,\notag\\
 B_\Gamma&\leq\|\lambda\|\,\|t\|_{\dot f^{s,q}_{p,F}}
 \leq\|\lambda\|\,|\Omega|^{1/p-1/q}.
\end{align}
The zero case is immediate. Thus
\begin{equation}\label{eq:carleson-reverse}
 \left\{\sum_{(j,k,R)\in\Gamma}
 |R|[(2^{-js_1}2^{-(j+k)s_2})^{-1}|d_{j,k,R}|]^{q'}\right\}^{1/q'}
 \leq\|\lambda\||\Omega|^{1/p-1/q}.
\end{equation}
Exhausting the labelled rectangles inside $\Omega$ and taking the
supremum over $\Omega$ proves
$\|d\|_{\dot c_F^{-s,p,q'}}\leq\|\lambda\|$. Finite sequences are
dense, so $\lambda(t)=[t,d]$ on the whole space, and the
representation is unique.
\end{proof}

\begin{proof}[Proof of Theorem~\ref{thm:carleson-dual}]
We first justify synthesis on the Carleson sequence space. A fixed
open concentric enlargement of $R$ has measure a fixed multiple of
$|R|$, so
\begin{equation}\label{eq:carleson-growth}
 |d_{j,k,R}|\leq C 2^{-js_1}2^{-(j+k)s_2}|R|^{1/p-1}
                         \|d\|_{\dot c_F^{-s,p,q'}}.
\end{equation}
The volume formula and \eqref{eq:testcoeff} show that
$T_\psi d$ converges absolutely against every test function and
defines a continuous element of $\mathcal S_F'$.

Put $P=S_\psi T_\psi$. By Theorem~\ref{thm:retraction}, $P$ is
bounded on $\dot f^{s,q}_{p,F}$. Its matrix is
\begin{equation}\label{eq:projection-matrix-expanded}
 (Pd)_{J,R}=\sum_{L,Q}|Q|\,
       (\psi_J*\psi_L)(z_R-z_Q)d_{L,Q}.
\end{equation}
For finite $t,d$, the real even kernels give
\begin{align}\label{eq:projection-transpose-expanded}
 [t,Pd]
 &=\sum_{J,R}\sum_{L,Q}|R||Q|\,
       t_{J,R}\overline{d_{L,Q}}\,
                (\psi_J*\psi_L)(z_R-z_Q)\notag\\
 &=\sum_{L,Q}|Q|(Pt)_{L,Q}\overline{d_{L,Q}}
 =[Pt,d].
\end{align} For
finite $d$, Lemma~\ref{lem:carleson-sequence} bounds the right side
by $C\|t\|_{\dot f^{s,q}_{p,F}}\|d\|_{\dot c_F^{-s,p,q'}}$. Taking the supremum over finitely supported $t$ and using
the converse sequence estimate gives
\begin{align}\label{eq:projection-dual-bound-expanded}
 \|Pd\|_{\dot c_F^{-s,p,q'}}
 &\leq C\sup_{\substack{t\ {\rm finite}\\
                    \|t\|_{\dot f^{s,q}_{p,F}}\leq1}}|[t,Pd]|\notag\\
 &=C\sup_{\substack{t\ {\rm finite}\\
                    \|t\|_{\dot f^{s,q}_{p,F}}\leq1}}|[Pt,d]|
 \leq C\|P\|\,\|d\|_{\dot c_F^{-s,p,q'}}.
\end{align} For general $d$, finite restrictions have no larger
Carleson norm and their syntheses converge in $\mathcal S_F'$, so
Fatou's lemma gives the same bound. Thus
$T_\psi:\dot c_F^{-s,p,q'}\to\mathrm{CMO}_F^{-s,p,q'}$ is bounded,
and $T_\psi S_\psi g=g$ by the exact formula.

If $\ell\in(\dot F^{s,q}_{p,F})^*$, apply sequence duality to
$\ell T_\psi$ to obtain $d\in\dot c_F^{-s,p,q'}$ and put
$g=T_\psi d$. The representation and its norm estimate are
\begin{align}\label{eq:carleson-representation-expanded}
 \ell(f)
 &=\ell(T_\psi S_\psi f)=[S_\psi f,d]
   =\langle f,T_\psi d\rangle=\langle f,g\rangle,
                     &&f\in\mathcal S_F,\notag\\
 \|g\|_{\mathrm{CMO}_{F}^{-s,p,q'}}
 &=\|Pd\|_{\dot c_F^{-s,p,q'}}
 \leq C\|d\|_{\dot c_F^{-s,p,q'}}
 \leq C\|\ell T_\psi\|
 \leq C\|\ell\|.
\end{align}
Density of $\mathcal S_F$ extends the representation to all $f$. Conversely,
\eqref{eq:carleson-pairing} shows that $[S_\psi f,S_\psi g]$ is a
bounded functional of $f$ and equals the distribution pairing on
test functions by the reproducing formula. This proves the duality
and the converse bound; uniqueness follows by testing against
$\mathcal S_F$. The continuous dual of a quasi-Banach space is
Banach, so the Carleson space is complete. Finally, for another
admissible pair and exact grid, the cross-system operator
$S_\varphi T_\psi$ and its adjoint $S_\psi T_\varphi$ are bounded
on the corresponding sequence spaces, by the analysis and synthesis
estimates; since $T_\psi S_\psi g=g$, this gives the independence
of the norm.
\end{proof}

\begin{proof}[Proof of Theorem~\ref{thm:duality}]
We first identify the dual of the coefficient space. To do this,
introduce the two step functions
\[
 U_J=\sum_R2^{-js_1}2^{-(j+k)s_2}|t_{J,R}|\mathbf1_R,
 \qquad
 V_J=\sum_R2^{js_1}2^{(j+k)s_2}|d_{J,R}|\mathbf1_R.
\]
Pointwise H\"older in the scales and then spatial H\"older give
\begin{align}\label{eq:sequence-holder}
 \sum_{J,R}|R|\,|t_{J,R}d_{J,R}|
 &=\int\sum_J U_J(z)V_J(z)\,dz\notag\\
 &\leq\int\left(\sum_J U_J^q\right)^{1/q}
             \left(\sum_J V_J^{q'}\right)^{1/q'}\,dz\notag\\
 &\leq\|t\|_{\dot f^{s,q}_{p,F}}\,
       \|d\|_{\dot f^{-s,q'}_{p',F}}.
\end{align}
This also proves absolute convergence. Conversely, let $\lambda$ be
a bounded functional on $\dot f^{s,q}_{p,F}$. Embed that space
isometrically into $L^p(\ell^q)$ by
$t\mapsto\{2^{-js_1}2^{-(j+k)s_2}\sum_R t_{j,k,R}\mathbf1_R\}$. Hahn--Banach and
the duality of this mixed Lebesgue space give
$H=\{H_{j,k}\}\in L^{p'}(\ell^{q'})$, with $\|H\|=\|\lambda\|$, such
that
$\lambda(t)=\sum_{j,k}\sum_R 2^{-js_1}2^{-(j+k)s_2} t_{j,k,R}\int_R\overline{H_{j,k}}$.
Define
\begin{equation}\label{eq:dual-coefficient}
 d_{j,k,R}=\frac{2^{-js_1}2^{-(j+k)s_2}}{|R|}\int_R H_{j,k}(z)\,dz.
\end{equation}
For $z\in R$,
\[
 2^{js_1}2^{(j+k)s_2}|d_{J,R}|
 \leq |R|^{-1}\int_R|H_J(u)|\,du
 \leq M_FH_J(z).
\]
Consequently Lemma~\ref{lem:maximal} gives
\begin{align}\label{eq:Banach-dual-average-expanded}
 \|d\|_{\dot f^{-s,q'}_{p',F}}
 &\leq\|\{M_FH_J\}_J\|_{L^{p'}(\ell^{q'})}\notag\\
 &\leq C\|H\|_{L^{p'}(\ell^{q'})}=C\|\lambda\|.
\end{align}

Let $\ell\in(\dot F^{s,q}_{p,F})^*$. Apply the coefficient duality
to $\ell T_\psi$ to obtain $d$ and set $g=T_\psi d$, which lies in
$\dot F^{-s,q'}_{p',F}$ by the synthesis theorem with parameters
$p',q',-s$. For $f\in\mathcal S_F$,
$\ell(f)=\ell(T_\psi S_\psi f)=[S_\psi f,d]=\langle f,g\rangle$.
Density extends this to the whole space and gives the norm estimate.
In the other direction, the exact formula gives
$\langle f,g\rangle=[S_\psi f,S_\psi g]$ for test functions, and
\eqref{eq:sequence-holder} extends this to a bounded pairing;
uniqueness follows by testing against $\mathcal S_F$.
\end{proof}

Write $\dot f^{s,\infty}_{p,F,0}$ for the closure of the finitely
supported sequences in $\dot f^{s,\infty}_{p,F}$.
The analysis and synthesis maps restrict to this closure and the
little function space: synthesis sends finite sequences into
$\mathcal S_F$, while the scale decay of test-function coefficients
and spatial truncation give $S_\psi\mathcal S_F\subset
\dot f^{s,\infty}_{p,F,0}$. Their boundedness then gives the
corresponding maps between the closures.

\begin{lemma}[Hardy-range endpoint sequences]\label{lem:Hardy-endpoint-sequences}
Under the hypotheses of Theorem~\ref{thm:Hardy-endpoint},
\begin{equation}\label{eq:Hardy-endpoint-sequences}
 (\dot f^{s,q}_{p,F})^*=\dot c_F^{-s,p,\infty}\quad(0<q\leq1),
 \qquad
 (\dot f^{s,\infty}_{p,F,0})^*=\dot c_F^{-s,p,1}.
\end{equation}
The first identity is isometric for the displayed sequence norms;
both pairings converge absolutely.
\end{lemma}
\begin{proof}
First suppose $q\leq1$ and let $t$ have finite support. For
nonnegative functions and $0<p\leq1$, the reverse Minkowski
inequality gives
\begin{equation}\label{eq:positive-reverse}
\begin{aligned}
 \sum_{J,R}2^{-js_1}2^{-(j+k)s_2}|t_{J,R}||R|^{1/p}
 &\leq\left\|\sum_{J,R}2^{-js_1}2^{-(j+k)s_2}|t_{J,R}|\mathbf1_R\right\|_{L^p}\\
 &\leq\left\|\{\sum_{J,R}[2^{-js_1}2^{-(j+k)s_2}|t_{J,R}|]^q\mathbf1_R\}^{1/q}\right\|_{L^p}.
\end{aligned}
\end{equation}
The first inequality follows from the concavity of $x\mapsto x^p$:
write $A_i=2^{-j_i s_1}2^{-(j_i+k_i)s_2}$,
$b_i=A_i|t_i||R_i|^{1/p}$, $B=\sum_i b_i$, and
$u_i=\mathbf1_{R_i}/|R_i|^{1/p}$. The case $B=0$ is immediate. Otherwise, concavity gives
\begin{align}\label{eq:reverse-Minkowski-concavity-expanded}
 \int\left(\sum_i\frac{b_i}{B}u_i\right)^p
 &\geq\sum_i\frac{b_i}{B}\int u_i^p
 =\sum_i\frac{b_i}{B}=1,\notag\\
 \left\|\sum_i b_i u_i\right\|_{L^p}&\geq B=\sum_i b_i.
\end{align} The second follows from
$\sum_i c_i\leq(\sum_i c_i^q)^{1/q}$. Consequently
\begin{equation}\label{eq:supremum-pairing}
 \sum_{J,R}|R|\,|t_{J,R}d_{J,R}|
 \leq\|d\|_{\dot c_F^{-s,p,\infty}}\|t\|_{\dot f^{s,q}_{p,F}}.
\end{equation}
Finite restrictions give absolute convergence for arbitrary $t$.
Conversely, a bounded functional $\lambda$ determines $d$ by
$\lambda(e_{J,R})=|R|\overline{d_{J,R}}$. Coordinate testing gives
\begin{align}\label{eq:endpoint-coordinate-test-expanded}
 |R|\,|d_{J,R}|
 &=|\lambda(e_{J,R})|
 \leq\|\lambda\|\,2^{-js_1}2^{-(j+k)s_2}|R|^{1/p},\notag\\
 |R|^{1-1/p}2^{js_1}2^{(j+k)s_2}|d_{J,R}|
 &\leq\|\lambda\|.
\end{align}
Taking the supremum gives
$\|d\|_{\dot c_F^{-s,p,\infty}}\leq\|\lambda\|$. Together with
\eqref{eq:supremum-pairing}, this proves equality of the dual norms.
Density of finite sequences extends the representation.

For the second identity, begin with finite $t$ and put
$G(z)=\sup_{J,R}2^{-js_1}2^{-(j+k)s_2}|t_{J,R}|\mathbf1_R(z)$ and
$\Omega_l=\{z:G(z)>2^l\}$. Group the nonzero coefficients in
$\mathcal B_l$ as in the proof of Lemma~\ref{lem:carleson-sequence}
and let $\Omega_l^*=\{M_F\mathbf1_{\Omega_l}>1/4\}$. Then
$|\Omega_l^*|\leq C|\Omega_l|$ and every label in $\mathcal B_l$
satisfies $2^{-js_1}2^{-(j+k)s_2}|t_{J,R}|\leq2^{l+1}$, so
\begin{align}
 \sum_{\mathcal B_l}|R|\,|t_{J,R}d_{J,R}|
 &\leq2^{l+1}\sum_{R\subset\Omega_l^*}|R|(2^{-js_1}2^{-(j+k)s_2})^{-1}|d_{J,R}|\notag\\
 &\leq C\|d\|_{\dot c_F^{-s,p,1}}2^l|\Omega_l|^{1/p}.
\end{align}
Since $p\leq1$, summing in $l$ and the layer-cake estimate give
\begin{equation}\label{eq:c-one-pairing}
 \sum_{J,R}|R|\,|t_{J,R}d_{J,R}|
 \leq C\|d\|_{\dot c_F^{-s,p,1}}
   (\sum_l2^{lp}|\Omega_l|)^{1/p}
 \leq C\|d\|_{\dot c_F^{-s,p,1}}\|t\|_{\dot f^{s,\infty}_{p,F}}.
\end{equation}
This extends to the little subspace. Conversely, let $\lambda$ be
bounded on the little sequence space and define $d$ by coordinate
testing. Fix an open $\Omega$ of finite positive measure and a
finite set $\Gamma$ of labels with $R\subset\Omega$. Take
$t_{J,R}=(2^{-js_1}2^{-(j+k)s_2})^{-1}d_{J,R}/|d_{J,R}|$ on $\Gamma$ when
$d_{J,R}\ne0$, and set all other entries to zero. For this finite sequence,
\begin{align}\label{eq:c-one-test-expanded}
 \sup_{J,R}2^{-js_1}2^{-(j+k)s_2}|t_{J,R}|\mathbf1_R
 &\leq\mathbf1_\Omega,\notag\\
 \|t\|_{\dot f^{s,\infty}_{p,F}}
 &\leq\|\mathbf1_\Omega\|_{L^p}=|\Omega|^{1/p}.
\end{align}
In particular $t$ belongs to the little sequence space, and
\begin{equation}\label{eq:c-one-reverse}
 \sum_{(J,R)\in\Gamma}|R|(2^{-js_1}2^{-(j+k)s_2})^{-1}|d_{J,R}|=\lambda(t)
       \leq\|\lambda\||\Omega|^{1/p}.
\end{equation}
Letting $\Gamma$ exhaust the labels in $\Omega$ and taking the
supremum over $\Omega$ proves the Carleson bound; finite sequences
are dense by the definition of
the little space, so the representation holds and is unique.
\end{proof}

\begin{proof}[Proof of Theorem~\ref{thm:Hardy-endpoint}]
Denote the relevant primal sequence space by $X$ and its sequence
dual by $C$; in both cases finite sequences are dense, analysis and
synthesis give a retraction, and $P=S_\psi T_\psi$ is a bounded
projection. Both endpoint Carleson norms imply
\begin{equation}\label{eq:endpoint-growth}
 |d_{J,R}|\leq C 2^{-js_1}2^{-(j+k)s_2}|R|^{1/p-1}\|d\|_C.
\end{equation}
For the $r=1$ norm this uses an open concentric enlargement of $R$;
for the supremum norm it is immediate. The volume formula and
\eqref{eq:testcoeff} then give absolute convergence of $T_\psi d$
against every test function, hence a continuous element of
$\mathcal S_F'$.

The transpose identity \eqref{eq:projection-transpose-expanded}
and Lemma~\ref{lem:Hardy-endpoint-sequences} give, for finite $d$,
\begin{align}\label{eq:endpoint-projection-expanded}
 \|Pd\|_C
 &\leq C\sup_{\substack{t\ {\rm finite}\\\|t\|_X\leq1}}|[t,Pd]|
 =C\sup_{\substack{t\ {\rm finite}\\\|t\|_X\leq1}}|[Pt,d]|\notag\\
 &\leq C\|P\|_{X\to X}\|d\|_C.
\end{align}
For a finite exhaustion $\Gamma_a$, coefficient restriction gives
$\|\mathbf1_{\Gamma_a}d\|_C\leq\|d\|_C$, and distributional
convergence gives $P(\mathbf1_{\Gamma_a}d)\to Pd$ coordinatewise.
Fatou's lemma, followed by the defining supremum, gives
\begin{equation}\label{eq:endpoint-Fatou-expanded}
 \|T_\psi d\|_{\mathrm{CMO}}
 =\|Pd\|_C
 \leq\liminf_{a\to\infty}\|P(\mathbf1_{\Gamma_a}d)\|_C
 \leq C\|d\|_C.
\end{equation}

Let $\ell$ be a bounded functional on the corresponding primal
function space. Then $\ell T_\psi$ has a representative $d\in C$;
set $g=T_\psi d$. For $f\in\mathcal S_F$,
$\ell(f)=[S_\psi f,d]=\langle f,g\rangle$, and density extends this
to every $f$ with $\|g\|_{\mathrm{CMO}}\leq C\|\ell\|$. Conversely,
the sequence lemma bounds $[S_\psi f,S_\psi g]$, which equals the
distribution pairing on test functions. This proves both duality
identities and the reverse bounds; uniqueness follows by testing
against $\mathcal S_F$. The continuous dual of a quasi-Banach space
is Banach, so the endpoint Carleson spaces are complete. For the
norm independence, the cross-system operator and its adjoint are
bounded on the endpoint Carleson sequence spaces; the positive
cross matrix has the same almost-orthogonality majorant as in the
matrix lemma, and the finite-sequence and Fatou argument gives the
bounds, as in the proof of Theorem~\ref{thm:carleson-dual}.
\end{proof}

\begin{lemma}[Small-$q$ Banach sequence duality]\label{lem:small-q-Banach-sequence}
Let $1<p<\infty$, $0<q\leq1$ and $s\in\mathbb R^2$. Then
\begin{equation}\label{eq:small-q-Banach-sequence}
 (\dot f^{s,q}_{p,F})^*=\dot f^{-s,\infty}_{p',F}
\end{equation}
with equivalent norms; the pairing converges absolutely.
\end{lemma}
\begin{proof}
For a sequence $d$ put $b_{J,R}=2^{js_1}2^{(j+k)s_2}|d_{J,R}|$ and
$D(z)=\sup_{J,R}b_{J,R}\mathbf1_R(z)$. The inclusion
$\ell^q\subset\ell^1$ and spatial H\"older give
\begin{equation}\label{eq:small-q-Banach-pairing}
 \sum_{J,R}|R|\,|t_{J,R}d_{J,R}|
 \leq\int D(z)\sum_{J,R}2^{-js_1}2^{-(j+k)s_2}|t_{J,R}|\mathbf1_R(z)\,dz
 \leq\|D\|_{L^{p'}}\|t\|_{\dot f^{s,q}_{p,F}}.
\end{equation}
Conversely, let $\lambda$ be a bounded functional and define $d$ by
$\lambda(e_{J,R})=|R|\overline{d_{J,R}}$. For a finite collection of
labels, order them so that $b_{J,R}$ decreases and retain a rectangle
$R$ when at least half of it lies outside the previously retained
rectangles; the new parts $E_{J,R}$ are pairwise disjoint and
$|E_{J,R}|\geq|R|/2$. Write $D_0=\sup_{\mathrm{retained}}b_{J,R}\mathbf1_R$.
On $E_{J,R}$, $D_0=b_{J,R}$, so
\begin{equation}\label{eq:greedy-norm}
 \|D_0\|_{L^{p'}}^{p'}=\sum_{\mathrm{retained}}b_{J,R}^{p'}|E_{J,R}|.
\end{equation}
If a rectangle was discarded, more than half of it is covered by
previously retained rectangles, so for $0<r<\min(1,p,p',q)$,
$b_{J,R}\leq2^{1/r}[M_F(D_0^r)(z)]^{1/r}$ for $z\in R$; the same
holds for a retained rectangle. Taking the supremum and using $p'/r>1$ gives
\begin{align}\label{eq:greedy-maximal-expanded}
 D_\Gamma(z)&\leq2^{1/r}[M_F(D_0^r)(z)]^{1/r},\notag\\
 \|D_\Gamma\|_{L^{p'}}^r
 &\leq2\|M_F(D_0^r)\|_{L^{p'/r}}
 \leq C\|D_0^r\|_{L^{p'/r}}
 =C\|D_0\|_{L^{p'}}^r.
\end{align}

On the retained labels choose
$t_{J,R}=2^{js_1}2^{(j+k)s_2}b_{J,R}^{p'-1}d_{J,R}/|d_{J,R}|$ and zero
elsewhere, with $t_{J,R}=0$ whenever $d_{J,R}=0$.
Since $R$ is a flag rectangle and $|E_{J,R}|\geq|R|/2$,
\[
 b_{J,R}^{p'-1}\mathbf1_R
 \leq2^{1/r}\bigg(M_F((b_{J,R}^{p'-1}\mathbf1_{E_{J,R}})^r)\bigg)^{1/r}.
\]
Lemma~\ref{lem:maximal}, applied to the labelled family $(J,R)$ with
exponents $p/r,q/r>1$, and disjointness of the sets $E_{J,R}$ give
\begin{align}
 \|t\|_{\dot f^{s,q}_{p,F}}
 &\leq C\|\{b_{J,R}^{p'-1}\mathbf1_{E_{J,R}}\}\|_{L^p(\ell^q)}
 =C\|D_0\|_{L^{p'}}^{p'/p},\\
 \lambda(t)&=\sum_{\mathrm{retained}}|R|b_{J,R}^{p'}
       \geq\|D_0\|_{L^{p'}}^{p'}.
\end{align}
If $\|D_0\|_{L^{p'}}>0$, the preceding inequalities give
\begin{align}\label{eq:greedy-division-expanded}
 \|D_0\|_{L^{p'}}^{p'}
 &\leq|\lambda(t)|
 \leq\|\lambda\|\,\|t\|_{\dot f^{s,q}_{p,F}}
 \leq C\|\lambda\|\,\|D_0\|_{L^{p'}}^{p'/p},\notag\\
 \|D_0\|_{L^{p'}}
 &=\|D_0\|_{L^{p'}}^{\,p'-p'/p}
 \leq C\|\lambda\|.
\end{align}
The case $\|D_0\|_{L^{p'}}=0$ is immediate. The exponents used here
are $(p'-1)p=p'$ and $p'-p'/p=1$. Increasing the finite
collections proves $d\in\dot f^{-s,\infty}_{p',F}$ with the required
bound; finite sequences are dense, so $\lambda(t)=[t,d]$ on the whole
space and the representation is unique.
\end{proof}

\begin{lemma}[Averaging endpoints]\label{lem:averaging}
Let $1<p<\infty$ and define
\begin{equation}\label{eq:Q}
 Q_J H(z)=\sum_{R\in\mathcal R_J^N}\frac{\mathbf1_R(z)}{|R|}\int_R H(u)\,du.
\end{equation}
Then $Q=\{Q_J\}_J$ is bounded on both $L^{p'}(\ell^1)$ and
$L^{p'}(\ell^\infty)$.
\end{lemma}
\begin{proof}
For the supremum norm, positivity and scalar boundedness give
\begin{align}\label{eq:averaging-sup-expanded}
 \|QG\|_{L^{p'}(\ell^\infty)}
 &\leq\left\|M_F\!\left(\sup_J|G_J|\right)\right\|_{L^{p'}}
 \leq C\|G\|_{L^{p'}(\ell^\infty)}.
\end{align} For the $\ell^1$ norm, take finitely many
coordinates and $h\geq0$ with $\|h\|_{L^p}\leq1$. The averaging
operators are self-adjoint, so
\begin{align}
 \int h(z)\sum_J|Q_JG_J(z)|\,dz
 &\leq\sum_J\int(Q_Jh)(u)|G_J(u)|\,du\notag\\
 &\leq\int M_Fh(u)\sum_J|G_J(u)|\,du
 \leq C\|\sum_J|G_J|\|_{L^{p'}}.
\end{align}
Scalar duality gives the $L^{p'}(\ell^1)$ bound; monotone
convergence covers general families.
\end{proof}

\begin{proof}[Proof of Theorem~\ref{thm:Banach-endpoint}]
Lemma~\ref{lem:small-q-Banach-sequence} establishes the first
coefficient identity. For the second, the little sequence space
$\dot f^{s,\infty}_{p,F,0}$ embeds isometrically into $L^p(c_0)$ by
multiplication of its step functions by $2^{-js_1}2^{-(j+k)s_2}$. Hahn--Banach
extends any functional to the mixed Lebesgue space, with
representing family in $L^{p'}(\ell^1)$. Define
\begin{equation}\label{eq:Banach-endpoint-coefficient}
 d_{J,R}=\frac{2^{-js_1}2^{-(j+k)s_2}}{|R|}\int_R H_J(z)\,dz.
\end{equation}
Then $\lambda(t)=[t,d]$ for finite $t$, and
$(2^{-js_1}2^{-(j+k)s_2})^{-1}d_{J,R}$ is the value of $Q_JH_J$ on $R$.
Lemma~\ref{lem:averaging} gives
\begin{align}\label{eq:little-Banach-average-expanded}
 \|d\|_{\dot f^{-s,1}_{p',F}}
 &=\|QH\|_{L^{p'}(\ell^1)}
 \leq C\|H\|_{L^{p'}(\ell^1)}
 =C\|\lambda\|.
\end{align}
Conversely, the absolute pairing is bounded by
\begin{align}\label{eq:little-Banach-holder-expanded}
 \sum_{J,R}|R|\,|t_{J,R}d_{J,R}|
 &\leq\int
 \left(\sup_{J,R}2^{-js_1}2^{-(j+k)s_2}|t_{J,R}|\mathbf1_R\right)
 \left(\sum_{J,R}2^{js_1}2^{(j+k)s_2}|d_{J,R}|\mathbf1_R\right)\,dz\notag\\
 &\leq\|t\|_{\dot f^{s,\infty}_{p,F}}\,
            \|d\|_{\dot f^{-s,1}_{p',F}}.
\end{align}
Density of finite sequences in the little space completes the
sequence duality.

For a bounded functional $\ell$ on either function space, apply the
sequence result to $\ell T_\psi$ and put $g=T_\psi d$; the retraction
with parameters $p',-s$ places $g$ in the right side of
\eqref{eq:q-one-dual} or \eqref{eq:little-Banach-dual}. For test
functions, $\ell(f)=[S_\psi f,d]=\langle f,g\rangle$, and density
extends the representation. In the converse direction,
$[S_\psi f,S_\psi g]$ is bounded by the coefficient pairing and
equals the distribution pairing on tests; uniqueness follows by
testing against $\mathcal S_F$. This proves both identities.
\end{proof}

\begin{proof}[Proof of Corollary~\ref{cor:small-q-Besov}]
On a nonzero kernel $k\leq j+1$, and
$|R|\asymp_N2^{-jn}2^{-km}$. With $1-1/p=-\delta$,
\begin{equation}\label{eq:Besov-weight}
 2^{js_1}2^{(j+k)s_2}|R|^{1-1/p}
 \asymp 2^{-jt_1}2^{-(j+k)t_2}.
\end{equation}
The two exponent identities are
\begin{align}\label{eq:Besov-exponents-expanded}
 -t_2&=s_2+m\delta,\notag\\
 -(t_1+t_2)
 &=s_1-(m-n)\delta+s_2+m\delta
 =s_1+s_2+n\delta.
\end{align}
Hence the weight in \eqref{eq:Besov-weight} equals, up to a fixed
grid constant,
$2^{j(s_1+s_2+n\delta)}2^{k(s_2+m\delta)}$. Thus the endpoint Carleson norm
\eqref{eq:c-infinity} is equivalent to the supremum of the samples
defining the Besov norm. To compare the sampled and continuous
suprema, insert the exact formula for $g$ into $\psi_J*g$; only
finitely many neighbouring scales contribute, and almost
orthogonality and $\sum_R|R|P_{a_L,b_L,M}(z-z_R)\leq C$ for $M>D$
give
\begin{equation}\label{eq:Besov-supremum}
 \sup_z 2^{-jt_1}2^{-(j+k)t_2}|\psi_J*g(z)|
 \leq C\sup_{L,R}2^{-ut_1}2^{-(u+v)t_2}|(\psi_L*g)(z_R)|.
\end{equation}
The reverse inequality is immediate. Theorem~\ref{thm:Hardy-endpoint}
then proves the identity and the independence of the pair.
\end{proof}

\section{The full supremum dual: proofs of Propositions~\ref{prop:full-supremum} and~\ref{prop:modulations}}
\label{sec:full-dual}

\begin{proof}[Proof of Proposition~\ref{prop:full-supremum}]
The sequence estimates above hold on the entire supremum sequence
space: in the Hardy range this is proved in
Lemma~\ref{lem:Hardy-endpoint-sequences}, and in the Banach spatial
range it is the $\ell^\infty$--$\ell^1$ and spatial H\"older
inequality. They give absolute convergence of
\eqref{eq:regular-pairing} and
$|\langle f,g\rangle_F|\leq C\|f\|_X\|g\|_Y$.

We verify pair independence without density in $X$. The positive
matrix obtained by replacing every entry of the cross-system operator
$S_\varphi T_\psi$ by its absolute value is bounded on the full
supremum sequence spaces: its entries have the almost-orthogonality
majorant used in Lemma~\ref{lem:matrix}, whose proof applies with
$q=\infty$. For two fixed grid resolutions the cross sum is
absolutely convergent,
\begin{equation}\label{eq:cross-sum}
 \sum_{I,Q}\sum_{J,R}|Q||R|\,|t_{J,R}|\,|d_{I,Q}|\,
   |(\varphi_I*\psi_J)(z_Q-z_R)|
 \leq C\|t\|_{\dot f^{s,\infty}_{p,F}(\psi)}\|d\|_{Y_{\rm seq}(\varphi)},
\end{equation}
where $Y_{\rm seq}$ is $\dot c_F^{-s,p,1}$ or
$\dot f^{-s,1}_{p',F}$. Fubini and the real even kernels give the
transpose identity, and since the exact formulas hold in
$\mathcal S_F'$,
\begin{equation}\label{eq:pair-independence}
 [S_\psi f,S_\psi g]_\psi
 =[S_\varphi T_\psi S_\psi f,S_\varphi g]_\varphi
 =[S_\varphi f,S_\varphi g]_\varphi.
\end{equation}
The dual membership of $g$ for both choices follows from the
little-space duality or from Plancherel--P\'olya. If $f\in\mathcal S_F$,
its sampled reproducing formula converges in the test-function
topology, and pairing with $g$ gives the distribution pairing; if
$g\in\mathcal S_F$, the weak formula for $f$ and the same coefficient
bounds apply.

Restrict $\ell$ to $X_0$ and let $g\in Y$ be its unique
little-space representative. Define
\[
 \ell_g(f)=\langle f,g\rangle_F,\qquad \ell_0=\ell-\ell_g.
\]
Then
\begin{align}\label{eq:regular-decomposition-estimate-expanded}
 \ell_0|_{X_0}&=0,\notag\\
 \|g\|_Y&\leq C\|\ell|_{X_0}\|\leq C\|\ell\|,\notag\\
 \|\ell_0\|&\leq\|\ell\|+\|\ell_g\|
              \leq\|\ell\|+C\|g\|_Y\leq C\|\ell\|,\\
 \|\ell\|&\leq\|\ell_g\|+\|\ell_0\|
              \leq C\|g\|_Y+\|\ell_0\|.\notag
\end{align}
If $\ell=\ell_g+\ell_0=\ell_{\widetilde g}+\widetilde\ell_0$, then
\[
 \langle f,g-\widetilde g\rangle_F=0\quad(f\in X_0)
 \quad\Longrightarrow\quad
 g=\widetilde g,\qquad \ell_0=\widetilde\ell_0.
\]
The implication is precisely the uniqueness in little-space
duality. This proves the decomposition and its norm equivalence.
\end{proof}

\begin{proof}[Proof of Proposition~\ref{prop:modulations}]
Choose an admissible pair for which the multipliers have an exact plateau. Take a nonnegative smooth function $\chi$ supported in
$[3/5,8/5]$ and positive on $(3/5,8/5)$, and in each relevant
dimension define the radial multiplier
\begin{equation}\label{eq:modulation-mother}
 \widehat\psi(\zeta)
 =\frac{\chi(|\zeta|)}
   {\{\sum_{h\in\mathbb Z}\chi(2^{-h}|\zeta|)^2\}^{1/2}}
 \quad(\zeta\ne0),
\end{equation}
extended by zero near the origin. On $4/5<|\zeta|<6/5$ only the
term $h=0$ occurs, so $\widehat\psi=1$ there and all other dilates
vanish. These multipliers define real even admissible mothers.

Let $e_x$ be the first coordinate unit vector in $\mathbb R^n$.
Choose a nonzero $\phi\in\mathcal S$ whose Fourier transform is
supported in a small neighbourhood of $(0,e_y)$; then
$\widehat\psi^{(2)}=1$ on this support and $\phi\in\mathcal S_F$ by
\eqref{eq:partialmoments}. Choose $j_0$ large and set
$j_\ell=j_0+10\ell$. For $(\xi,\eta)\in\operatorname{supp}\widehat\phi$
and every $\ell$,
\begin{equation}\label{eq:modulation-support}
 \frac45<|2^{-j_\ell}(\xi+2^{j_\ell}e_x,\eta)|<\frac65.
\end{equation}
Thus the Fourier support of $e^{i2^{j_\ell}x_1}\phi$ lies in the
region where the $(j_\ell,0)$ multiplier equals one and all other
multipliers vanish. For $b\in\ell^\infty(\mathbb N)$ define
\begin{equation}\label{eq:modulation-example}
 Eb=\sum_{\ell=1}^\infty b_\ell 2^{j_\ell(s_1+s_2)}
       e^{i2^{j_\ell}x_1}\phi.
\end{equation}
For $h\in\mathcal S_F$, integration by parts in $x_1$ gives,
for every positive integer $M$,
\begin{align}\label{eq:modulation-convergence-expanded}
 |\langle e^{i2^{j_\ell}x_1}\phi,h\rangle|
 &\leq 2^{-Mj_\ell}
          \|\partial_{x_1}^M(\phi\overline h)\|_{L^1}
 \leq C_M2^{-Mj_\ell}\mathfrak p_M(h),\notag\\
 \sum_{\ell\geq1}|b_\ell|2^{j_\ell(s_1+s_2)}
       |\langle e^{i2^{j_\ell}x_1}\phi,h\rangle|
 &\leq C_M\|b\|_{\ell^\infty}\mathfrak p_M(h)
       \sum_{\ell\geq1}2^{-j_\ell(M-|s_1+s_2|)}<\infty.
\end{align}
Choose $M>|s_1+s_2|+1$. This defines $Eb$ continuously in
$\mathcal S_F'$. Frequency localisation gives
\begin{align}
 \psi_{j_\ell,0}*Eb
 &=b_\ell 2^{j_\ell(s_1+s_2)}e^{i2^{j_\ell}x_1}\phi,\\
 \psi_{j,k}*Eb&=0\quad\text{if }(j,k)\notin\{(j_\ell,0)\}.
\end{align}
Hence, for this admissible norm,
\begin{equation}\label{eq:modulation-norm-exact}
 \|Eb\|_X=\|\phi\|_{L^p}\sup_\ell|b_\ell|,\qquad
 \omega^F_{Eb,\infty}(K)=\|\phi\|_{L^p}\sup_{j_\ell>K}|b_\ell|.
\end{equation}
The lower distance bound follows from Theorem~\ref{thm:little}; for
the upper bound, truncate \eqref{eq:modulation-example} after
finitely many terms. Each truncation is in $\mathcal S_F$ and the
remainder has norm $\|\phi\|_{L^p}\sup_{\ell>L}|b_\ell|$, so the
distance is $\|\phi\|_{L^p}\limsup_\ell|b_\ell|$. Independence of
the admissible pair gives
\eqref{eq:modulation-norm}--\eqref{eq:modulation-distance}. The
constant sequence $b_\ell=1$ gives $Eb\notin X_0$, and the embedding
of $\ell^\infty$ proves nonseparability.

If $p\geq1$, Theorem~\ref{thm:retraction} makes $X$ a Banach space.
Let $\mathbf1=(1,1,\ldots)$ and $x=E\mathbf1+X_0$. The distance
formula gives
\[
 \|x\|_{X/X_0}=\operatorname{dist}_X(E\mathbf1,X_0)>0.
\]
The functional $L_0(cx)=c$ on the line spanned by $x$ is bounded,
with $\|L_0\|=\|x\|_{X/X_0}^{-1}$. Hahn--Banach extends it to
$L\in(X/X_0)^*$. If $\pi_0:X\to X/X_0$ is the quotient map, then
\begin{equation}\label{eq:annihilator-explicit-expanded}
 \ell_0=L\pi_0\in X_0^\perp,\qquad
 \ell_0(E\mathbf1)=1,\qquad
 \|\ell_0\|\leq\|L\|.
\end{equation}
Thus the annihilator is nonzero. A regular functional that vanishes on
$X_0$ must have $g=0$ by the little-space duality, so the annihilator
functionals are additional to the regular ones.
\end{proof}

\section{Interpolation: proofs of Theorems~\ref{thm:complex-interpolation} and~\ref{thm:real-interpolation}}
\label{sec:interpolation}

For a coordinatewise measurable family $H=\{H_J\}_{J\in\mathbb Z^2}$
put
\begin{equation}\label{eq:E}
 \|H\|_{E_p^{s,q}}
 =\left\|\{\sum_J[2^{-js_1}2^{-(j+k)s_2}|H_J|]^q\}^{1/q}\right\|_{L^p},
\end{equation}
with the inner expression a supremum when $q=\infty$.

\begin{lemma}[Common continuous retraction]\label{lem:continuous-retraction}
Let $1<p<\infty$, $1\leq q\leq\infty$ and $s\in\mathbb R^2$. The maps
\begin{equation}\label{eq:continuous-maps}
 \mathcal A f=\{\psi_J*f\}_J,\qquad
 \mathcal B H=\sum_J\psi_J*H_J
\end{equation}
satisfy
$\mathcal A:\dot F^{s,q}_{p,F}\to E_p^{s,q}$ and
$\mathcal B:E_p^{s,q}\to\dot F^{s,q}_{p,F}$, with
$\mathcal B\mathcal A=I$; analysis is isometric for the fixed pair
and synthesis is bounded.
\end{lemma}
\begin{proof}
Analysis has the asserted norm by definition. To define synthesis,
let $q'$ be the conjugate exponent, with $q'=\infty$ for $q=1$ and
$q'=1$ for $q=\infty$. For $\phi\in\mathcal S_F$, H\"older in the
indices and then in the spatial variables gives
\begin{equation}\label{eq:continuous-pairing}
\begin{aligned}
 \sum_J\int|H_J(u)|\,|(\psi_J*\phi)(u)|\,du
 &\leq\|H\|_{E_p^{s,q}}
 \|\{(2^{-js_1}2^{-(j+k)s_2})^{-1}|\psi_J*\phi|\}_J\|_{L^{p'}(\ell^{q'})}\\
 &\leq C\|H\|_{E_p^{s,q}}\mathfrak p(\phi),
\end{aligned}
\end{equation}
where the last inequality follows from \eqref{eq:testcoeff}. This
defines $\mathcal B H$ continuously in $\mathcal S_F'$.

For these annular kernels, $\psi_J*\psi_L=0$ unless
$|j-l_1|+|k-l_2|\leq c_0$; write $L\sim J$ for these finitely many
neighbours. Almost orthogonality gives, for $T>D$,
\begin{equation}\label{eq:continuous-kernel}
 |\psi_J*\psi_L(x,y)|
 \leq C_T a_L^{-n}b_L^{-m}
   (1+\max\{|x|/a_L,|y|/b_L\})^{-T}\qquad(L\sim J).
\end{equation}
For the positive kernel on the right, write
$\rho_L(x,y)=\max(|x|/a_L,|y|/b_L)$. Decomposition into
$\rho_L$-shells gives, for $g\geq0$,
\begin{align}\label{eq:continuous-shells-expanded}
 &\int a_L^{-n}b_L^{-m}(1+\rho_L(v))^{-T}g(z-v)\,dv\notag\\
 &\qquad\leq C\sum_{h\geq0}2^{-h(T-D)}
 \frac{1}{|E_{L,h}|}\int_{E_{L,h}}g(z-v)\,dv
 \leq C M_Fg(z),
\end{align}
where $E_{L,h}$ has side lengths comparable to $2^ha_L,2^hb_L$.
The sum converges because $T>D$. Also, for $L=(u,v)\sim J$,
\[
 \frac{2^{-js_1}2^{-(j+k)s_2}}{2^{-us_1}2^{-(u+v)s_2}}
 \leq2^{c_0(|s_1+s_2|+|s_2|)}.
\]
Consequently
\begin{equation}\label{eq:continuous-majorant}
 2^{-js_1}2^{-(j+k)s_2}|\psi_J*\mathcal B H|
 \leq C_s\sum_{L\sim J}M_F(2^{-us_1}2^{-(u+v)s_2}|H_L|).
\end{equation}
For $1<q<\infty$, the number of neighbours is uniformly bounded,
so Lemma~\ref{lem:maximal} gives
\begin{align}\label{eq:continuous-synthesis-expanded}
 \|\mathcal B H\|_{\dot F^{s,q}_{p,F}}
 &\leq C_s\|\{M_F(2^{-js_1}2^{-(j+k)s_2}|H_J|)\}_J
                                  \|_{L^p(\ell^q)}\notag\\
 &\leq C\|H\|_{E_p^{s,q}}.
\end{align}
For $q=\infty$, the corresponding estimate is
\[
 \sup_J2^{-js_1}2^{-(j+k)s_2}|\psi_J*\mathcal B H|
 \leq C_s M_F\!\left(\sup_L2^{-us_1}2^{-(u+v)s_2}|H_L|\right),
\]
followed by scalar boundedness on $L^p$. For $q=1$, scalar duality gives the bound: a
nonnegative kernel $K_L$ majorising the finitely many convolution
kernels satisfies $K_L*h\leq C M_Fh$, so for $h\geq0$,
$\|h\|_{L^{p'}}\leq1$,
\begin{equation}\label{eq:q-one-synthesis}
 \int h\sum_J 2^{-js_1}2^{-(j+k)s_2}|\psi_J*\mathcal B H|
 \leq C\sum_L\int(K_L*h)2^{-us_1}2^{-(u+v)s_2}|H_L|
 \leq C\|H\|_{E_p^{s,1}}.
\end{equation}
The identity $\mathcal B\mathcal A=I$ is
\eqref{eq:continuousCRF}. The maps do not depend on the norm
parameters, and the same test pairings give the compatibility on
sums of the spaces.
\end{proof}

For a compatible couple $(X_0,X_1)$ we use
\begin{equation}\label{eq:Kdef}
 K(t,h;X_0,X_1)=\inf_{h=h_0+h_1}(\|h_0\|_{X_0}+t\|h_1\|_{X_1}),
\end{equation}
and write $(X_0,X_1)_{\theta,r}$ for the real interpolation space.
Complex interpolation is in the sense of Calder\'on
\cite{Calderon64Interpolation}.

\begin{proposition}[$K$-characterisation]\label{prop:K}
Let $1<p_0,p_1<\infty$, $1\leq q_0,q_1\leq\infty$ and
$s^0,s^1\in\mathbb R^2$. Set $X_i=\dot F^{s^i,q_i}_{p_i,F}$ and
$E_i=E_{p_i}^{s^i,q_i}$. Then
\begin{equation}\label{eq:K-characterisation}
 K(t,f;X_0,X_1)\asymp K(t,\mathcal A f;E_0,E_1).
\end{equation}
\end{proposition}
\begin{proof}
For any decomposition $f=f_0+f_1$, isometric analysis gives
\begin{align}\label{eq:K-analysis-expanded}
 K(t,\mathcal A f;E_0,E_1)
 &\leq\|\mathcal A f_0\|_{E_0}
               +t\|\mathcal A f_1\|_{E_1}\notag\\
 &=\|f_0\|_{X_0}+t\|f_1\|_{X_1}.
\end{align}
Taking the infimum gives $K(t,\mathcal A f;E_0,E_1)\leq K(t,f;X_0,X_1)$.
Conversely, if $\mathcal A f=H_0+H_1$, then
\begin{align}\label{eq:K-synthesis-expanded}
 f&=\mathcal B H_0+\mathcal B H_1,\notag\\
 K(t,f;X_0,X_1)
 &\leq\|\mathcal B H_0\|_{X_0}
          +t\|\mathcal B H_1\|_{X_1}\notag\\
 &\leq\max_{i=0,1}\|\mathcal B\|_{E_i\to X_i}
            \bigl(\|H_0\|_{E_0}+t\|H_1\|_{E_1}\bigr).
\end{align}
The infimum over $H_0,H_1$ proves the reverse bound.
\end{proof}

\begin{lemma}[Mixed complex interpolation]\label{lem:mixed-complex}
Let $1<p_0,p_1,q_0,q_1<\infty$ and $0<\theta<1$, with the parameters
in \eqref{eq:interpolated-parameters}. Then
\begin{equation}\label{eq:mixed-complex}
 [E_{p_0}^{s^0,q_0},E_{p_1}^{s^1,q_1}]_\theta=E_p^{s,q}
\end{equation}
with equal norms.
\end{lemma}
\begin{proof}
Set
\[
 W_J^0=2^{-js_1^0}2^{-(j+k)s_2^0},\qquad
 W_J^1=2^{-js_1^1}2^{-(j+k)s_2^1},
\]
so the weight at the interpolated smoothness is
\[
 W_J=(W_J^0)^{1-\theta}(W_J^1)^\theta
     =2^{-js_1}2^{-(j+k)s_2}.
\]
First take a nonzero simple family $H$ with finitely many coordinates
and finite spatial support. By homogeneity, assume that
$\|H\|_{E_p^{s,q}}=1$, and put
$y_J=W_J H_J$, $B=\{\sum_J|y_J|^q\}^{1/q}$ and
$r_J=|y_J|/B$ on $\{B>0\}$, setting $r_J=0$ on $\{B=0\}$.
Then $\int B^p=1$ and $\sum_J r_J^q=1$ on $\{B>0\}$.
For $0\leq\Re\zeta\leq1$ define
\begin{equation}\label{eq:analytic-family}
 \frac1{p(\zeta)}=\frac{1-\zeta}{p_0}+\frac\zeta{p_1},\quad
 \frac1{q(\zeta)}=\frac{1-\zeta}{q_0}+\frac\zeta{q_1},
\end{equation}
and the endpoint weight at $\zeta$ is
\begin{equation}\label{eq:analytic-weight}
 W_J(\zeta)=(W_J^0)^{1-\zeta}(W_J^1)^\zeta,
\end{equation}
and
\begin{equation}\label{eq:analytic-entry}
 F_J(\zeta)=W_J(\zeta)^{-1}B^{p/p(\zeta)}
       r_J^{q/q(\zeta)}H_J/|H_J|.
\end{equation}
Set $F_J(\zeta)=0$ when $H_J=0$. The family is analytic and
$F(\theta)=H$. On the boundary $\nu=0,1$, all imaginary powers have
modulus one, and
\begin{align}\label{eq:analytic-boundary-expanded}
 |W_J^\nu F_J(\nu+iv)|&=B^{p/p_\nu}r_J^{q/q_\nu},\notag\\
 \sum_J|W_J^\nu F_J(\nu+iv)|^{q_\nu}
 &=B^{pq_\nu/p_\nu}\sum_Jr_J^q=B^{pq_\nu/p_\nu},\notag\\
 \|F(\nu+iv)\|_{E_{p_\nu}^{s^\nu,q_\nu}}^{p_\nu}
 &=\int B^p=1.
\end{align}
Multiplication by $\exp(\delta(\zeta-\theta)^2)$ imposes decay as
$|\Im\zeta|\to\infty$, while
\[
 \sup_v\left|\exp(\delta(\nu+iv-\theta)^2)\right|
 \leq\exp\!\left(\delta\max(\theta^2,(1-\theta)^2)\right)
 \longrightarrow1 .
\]
This proves the upper interpolation estimate.

For the converse, the bilinear dual mixed norm is
\[
 \|G\|_{(E_p^{s,q})^\dagger}
 =\|\{W_J^{-1}G_J\}_J\|_{L^{p'}(\ell^{q'})}.
\]
Given a simple finite $G$ of dual norm one, the same construction
with $p_i',q_i'$ and reciprocal weights gives an analytic family
$\mathcal G$ through $G$, with endpoint dual norms one up to an
arbitrarily small factor. For an admissible family $\mathcal F$
through $H$, set
\[
 \Phi(\zeta)=\int\sum_J\mathcal F_J(\zeta,z)
                              \mathcal G_J(\zeta,z)\,dz .
\]
Boundary H\"older and the three-lines theorem yield
\begin{align}\label{eq:complex-reverse-expanded}
 |\Phi(\nu+iv)|
 &\leq\|\mathcal F(\nu+iv)\|_{E_{p_\nu}^{s^\nu,q_\nu}}\,
       \|\mathcal G(\nu+iv)\|_{(E_{p_\nu}^{s^\nu,q_\nu})^\dagger},\notag\\
 \left|\int\sum_JH_JG_J\right|
 &\leq\|\mathcal F\|_{\mathcal F(E_0,E_1)}.
\end{align}
Letting the decay factors tend to one, taking the supremum over
$G$, and then the infimum over $\mathcal F$ gives the reverse bound.
Simple finite families are dense because all four endpoint
exponents are finite and the spatial measure is sigma-finite.
\end{proof}

\begin{proof}[Proof of Theorem~\ref{thm:complex-interpolation}]
Write $X_i=\dot F^{s^i,q_i}_{p_i,F}$, $E_i=E_{p_i}^{s^i,q_i}$,
and $X=\dot F^{s,q}_{p,F}$. Interpolating the endpoint maps gives
\[
 \mathcal A:[X_0,X_1]_\theta\longrightarrow[E_0,E_1]_\theta,
 \qquad
 \mathcal B:[E_0,E_1]_\theta\longrightarrow[X_0,X_1]_\theta .
\]
Lemma~\ref{lem:mixed-complex} and
$\mathcal B\mathcal A=I$ give the two estimates
\begin{align}\label{eq:complex-retraction-expanded}
 \|f\|_X
 &=\|\mathcal A f\|_{E_p^{s,q}}
 =\|\mathcal A f\|_{[E_0,E_1]_\theta}
 \leq C\|f\|_{[X_0,X_1]_\theta},\notag\\
 \|f\|_{[X_0,X_1]_\theta}
 &=\|\mathcal B\mathcal A f\|_{[X_0,X_1]_\theta}
 \leq C\|\mathcal A f\|_{[E_0,E_1]_\theta}
 =C\|f\|_X.
\end{align}
These also prove both set inclusions.
\end{proof}

\begin{proof}[Proof of Theorem~\ref{thm:real-interpolation}]
For an arbitrary family $H$, write
$G=\|\{2^{-js_1}2^{-(j+k)s_2} H_J\}_J\|_{\ell^q}$. The mixed-space $K$-functional
equals the scalar one:
\begin{equation}\label{eq:mixed-K}
 K(t,H;E_{p_0}^{s,q},E_{p_1}^{s,q})
 =K(t,G;L^{p_0},L^{p_1}).
\end{equation}
For a positive splitting $G=g_0+g_1$, set
$H_J^i=(g_i/G)H_J$ on $\{G>0\}$ and zero elsewhere. Then
\begin{align}\label{eq:radial-splitting-expanded}
 H&=H^0+H^1,\notag\\
 \|\{2^{-js_1}2^{-(j+k)s_2}H_J^i(z)\}_J\|_{\ell^q}
 &=\frac{g_i(z)}{G(z)}
       \|\{2^{-js_1}2^{-(j+k)s_2}H_J(z)\}_J\|_{\ell^q}
 =g_i(z),\notag\\
 \|H^i\|_{E_{p_i}^{s,q}}&=\|g_i\|_{L^{p_i}}.
\end{align}
This proves the mixed $K$-functional is at most the scalar one.
Conversely, a splitting $H=H^0+H^1$ gives
$G\leq G_0+G_1$, where
$G_i=\|\{2^{-js_1}2^{-(j+k)s_2}H_J^i\}_J\|_{\ell^q}$.
Define
\begin{equation}\label{eq:scalar-redistribution-expanded}
 g_i=\frac{GG_i}{G_0+G_1}\quad(G_0+G_1>0),\qquad g_i=0
 \quad(G_0+G_1=0).
\end{equation}
Then $G=g_0+g_1$, $0\leq g_i\leq G_i$, and
$\|g_i\|_{L^{p_i}}\leq\|H^i\|_{E_{p_i}^{s,q}}$.
Taking the two infima proves \eqref{eq:mixed-K}.

The scalar real interpolation identity
\begin{equation}\label{eq:scalar-Lp}
 (L^{p_0},L^{p_1})_{\theta,p}=L^p
\end{equation}
follows from the classical rearrangement formula; see
\cite[Chapter~5]{BL76}. More explicitly, put
$\sigma=(1/p_0-1/p_1)^{-1}$. If $G^*$ denotes the decreasing
rearrangement, then
\begin{align}\label{eq:scalar-K-rearrangement-expanded}
 K(t,G;L^{p_0},L^{p_1})
 &\asymp\left(\int_0^{t^\sigma}G^*(u)^{p_0}\,du\right)^{1/p_0}
       +t\left(\int_{t^\sigma}^{\infty}G^*(u)^{p_1}\,du\right)^{1/p_1}.
\end{align}
Hardy's inequalities applied to the two integrals give
\begin{equation}\label{eq:real-norm-expanded}
 \left(\int_0^\infty
       [t^{-\theta}K(t,G;L^{p_0},L^{p_1})]^p\frac{dt}{t}\right)^{1/p}
 \asymp\left(\int_0^\infty G^*(u)^p\,du\right)^{1/p}
 =\|G\|_{L^p}.
\end{equation}
Thus Proposition~\ref{prop:K} and \eqref{eq:mixed-K} yield
\begin{align}\label{eq:real-interpolation-chain-expanded}
 \|f\|_{(\dot F^{s,q}_{p_0,F},\dot F^{s,q}_{p_1,F})_{\theta,p}}
 &\asymp\left(\int_0^\infty
       [t^{-\theta}K(t,G_f;L^{p_0},L^{p_1})]^p\frac{dt}{t}\right)^{1/p}
 \asymp\|G_f\|_{L^p}
 =\|f\|_{\dot F^{s,q}_{p,F}}.
\end{align}
This proves both \eqref{eq:scalar-K} and the asserted norm identity. For the reverse inclusion, if
$G_f\in L^p$, the scalar splitting and \eqref{eq:mixed-K} place
$\mathcal A f$ in the mixed endpoint sum, and synthesis puts $f$ in
the function-space endpoint sum without a prior compatibility
assumption.
\end{proof}

\section{Lifting: proof of Theorem~\ref{thm:lifting}}
\label{sec:lifting}

Multiplication by $m_\alpha$ and $m_{-\alpha}$ preserves
$\widehat{\mathcal S_F}$ continuously: the derivatives of these
functions have at most polynomial growth at infinity and fixed
power singularities at $\eta=0$, and the flatness defining
$\widehat{\mathcal S_F}$ absorbs the latter. Thus $\mathcal D^\alpha$
is defined on $\mathcal S_F$ and on $\mathcal S_F'$ by transposition.

\begin{lemma}[Lifting kernel]\label{lem:lifting-kernel}
Write $J=(j,k)$ and $L=(l_1,l_2)$. The kernel
$\psi_J*\mathcal D^\alpha\psi_L$ vanishes unless
$|j-l_1|+|k-l_2|\leq c_0$, with $c_0$ independent of the indices.
For every $M>0$,
\begin{equation}\label{eq:lifting-kernel}
 |\psi_J*\mathcal D^\alpha\psi_L(z)|
 \leq C_{\alpha,M}2^{j(\alpha_1+\alpha_2)+k\alpha_2}
   P_{a_j,b_{j,k},M}(z)\qquad(L\sim J).
\end{equation}
\end{lemma}
\begin{proof}
On the nonempty Fourier support $|(\xi,\eta)|\asymp2^j$,
$|\eta|\asymp2^k$ and $k\leq j+c_0$. First take $k\leq j$ and
rescale $(\xi,\eta)=(2^j\zeta,2^k\nu)$. The normalised multiplier is
\begin{align}\label{eq:lift-normalised-symbol-expanded}
 a_{J,L}(\zeta,\nu)
 &=\widehat\psi_J(2^j\zeta,2^k\nu)
   \widehat\psi_L(2^j\zeta,2^k\nu)\notag\\
 &\qquad{}\times|(\zeta,2^{k-j}\nu)|^{\alpha_1+\alpha_2}
                       |\nu|^{\alpha_2}.
\end{align}
On its support,
\[
 |\zeta|\leq C,\qquad c\leq|\nu|\leq C,\qquad
 c\leq|(\zeta,2^{k-j}\nu)|\leq C,\qquad 0<2^{k-j}\leq1.
\]
Consequently, for every pair of multi-indices,
\begin{equation}\label{eq:lift-symbol-derivatives-expanded}
 \sup_{J,L:\,L\sim J}
 \|\partial_\zeta^\beta\partial_\nu^\gamma a_{J,L}\|_\infty
 \leq C_{\alpha,\beta,\gamma}.
\end{equation}
This includes negative powers because both radial factors are
bounded away from zero. Fourier inversion and integration by parts
give
\begin{align}\label{eq:lift-kernel-decay-expanded}
 |\psi_J*\mathcal D^\alpha\psi_L(x,y)|
 &\leq C_{\alpha,M}
 2^{j(\alpha_1+\alpha_2)+k\alpha_2}2^{jn+km}\notag\\
 &\qquad{}\times(1+2^j|x|)^{-n-M}(1+2^k|y|)^{-m-M}.
\end{align}
This is \eqref{eq:lifting-kernel} when $k\leq j$. If $j<k\leq j+c_0$, the two lengths are
comparable and the isotropic rescaling gives the same estimate.
\end{proof}

\begin{proof}[Proof of Theorem~\ref{thm:lifting}]
Apply $\mathcal D^\alpha$ to the exact discrete formula for $f$ and
convolve with $\psi_J$. Lemma~\ref{lem:lifting-kernel} gives
\begin{equation}\label{eq:lifting-apply}
 |\psi_J*\mathcal D^\alpha f(z)|
 \leq C\sum_{L\sim J}\sum_{R\in\mathcal R_L^N}|R|\,
   2^{j(\alpha_1+\alpha_2)+k\alpha_2}
   P_{a_j,b_{j,k},M}(z-z_R)|(\psi_L*f)(z_R)|.
\end{equation}
The smoothness factors satisfy
\begin{equation}\label{eq:lifting-cancel}
 2^{-j(s_1+\alpha_1)}2^{-(j+k)(s_2+\alpha_2)}
 2^{j(\alpha_1+\alpha_2)+k\alpha_2}
 =2^{-js_1}2^{-(j+k)s_2}.
\end{equation}
Only finitely many neighbouring scales occur. Thus
Lemma~\ref{lem:matrix}, \eqref{eq:lifting-cancel}, and analysis give
\begin{align}\label{eq:lift-norm-expanded}
 \|\mathcal D^\alpha f\|_{\dot F^{s+\alpha,q}_{p,F}}
 &\leq C\|S_\psi f\|_{\dot f^{s,q}_{p,F}}
 \leq C\|f\|_{\dot F^{s,q}_{p,F}}.
\end{align}
Finite restrictions and Fatou's lemma justify this for arbitrary
coefficients. Since $m_{-\alpha}m_\alpha=1$ on $\eta\ne0$,
\begin{align}\label{eq:lift-inverse-expanded}
 \mathcal D^{-\alpha}\mathcal D^\alpha f
 &=f=\mathcal D^\alpha\mathcal D^{-\alpha}f,\notag\\
 \|f\|_{\dot F^{s,q}_{p,F}}
 &=\|\mathcal D^{-\alpha}\mathcal D^\alpha f\|_{\dot F^{s,q}_{p,F}}
 \leq C\|\mathcal D^\alpha f\|_{\dot F^{s+\alpha,q}_{p,F}}.
\end{align}
Finally, if $f_a\in\mathcal S_F$ and $f_a\to f$ in
$\dot F^{s,\infty}_{p,F}$, then
\begin{equation}\label{eq:lift-little-expanded}
 \mathcal D^\alpha f_a\in\mathcal S_F,\qquad
 \|\mathcal D^\alpha f-\mathcal D^\alpha f_a\|_{\dot F^{s+\alpha,\infty}_{p,F}}
 \leq C\|f-f_a\|_{\dot F^{s,\infty}_{p,F}}\longrightarrow0.
\end{equation}
Applying the same argument to $-\alpha$ proves the little-space
isomorphism.
\end{proof}

\section{Embeddings: proof of Proposition~\ref{prop:embeddings} and Corollary~\ref{cor:Hardy-Besov}}
\label{sec:embeddings}

Use the same admissible pair in both norms and put
$G_J=2^{-js_1}2^{-(j+k)s_2}|\psi_{j,k}*f|$.

\begin{proof}[Proof of Proposition~\ref{prop:embeddings}]
All assertions follow from mixed-norm inequalities for nonnegative
functions $G_J$. First let $q\leq p$. Minkowski's inequality in
$L^{p/q}$ gives
\begin{equation}\label{eq:embedding-first}
 \|\{G_J\}\|_{L^p(\ell^q)}^q
 =\|\sum_J G_J^q\|_{L^{p/q}}
 \leq\sum_J\|G_J\|_{L^p}^q,
\end{equation}
which proves the first inclusion. The other inclusion for
$q\leq p$ follows from
\begin{align}\label{eq:embedding-q-small-expanded}
 \sum_J\|G_J\|_{L^p}^p
 &=\int\sum_JG_J(z)^p\,dz
 \leq\int\left(\sum_JG_J(z)^q\right)^{p/q}\,dz
 =\|\{G_J\}\|_{L^p(\ell^q)}^p.
\end{align}
If $p\leq q<\infty$, the first inclusion follows from
\begin{equation}\label{eq:embedding-p-small-expanded}
 \|\{G_J\}\|_{L^p(\ell^q)}^p
 =\int\left(\sum_JG_J(z)^q\right)^{p/q}\,dz
 \leq\int\sum_JG_J(z)^p\,dz
 =\sum_J\|G_J\|_{L^p}^p.
\end{equation}
Minkowski's integral inequality in $\ell^{q/p}$ gives the other
inclusion:
\begin{equation}\label{eq:embedding-second}
 \|\{\|G_J\|_{L^p}\}\|_{\ell^q}^{p}
 \leq\int\|\{G_J(z)^p\}\|_{\ell^{q/p}}dz
 =\|\{G_J\}\|_{L^p(\ell^q)}^p.
\end{equation}
For $q=\infty$ the complete comparison is
\begin{equation}\label{eq:embedding-sup-expanded}
 \sup_J\|G_J\|_{L^p}
 \leq\left\|\sup_JG_J\right\|_{L^p}
 \leq\left(\sum_J\|G_J\|_{L^p}^p\right)^{1/p}.
\end{equation}
Every use of Minkowski has exponent at least one, so the argument
covers all $p>0$. When $p=q$,
\[
 \|\{G_J\}\|_{L^p(\ell^p)}^p
 =\int\sum_JG_J^p
 =\sum_J\int G_J^p
 =\|\{\|G_J\|_{L^p}\}_J\|_{\ell^p}^p.
\]
Finally, for $q_0\leq q_1$,
\[
 \|\{G_J\}\|_{L^p(\ell^{q_1})}
 \leq\|\{G_J\}\|_{L^p(\ell^{q_0})},\qquad
 \|\{\|G_J\|_{L^p}\}_J\|_{\ell^{q_1}}
 \leq\|\{\|G_J\|_{L^p}\}_J\|_{\ell^{q_0}}.
\]
These are precisely \eqref{eq:q-monotonicity}.
\end{proof}

\begin{proof}[Proof of Corollary~\ref{cor:Hardy-Besov}]
The flag square-function theory of Han and Lu \cite{HL08,HL10}
identifies $\dot F^{0,2}_{p,F}$ with $H_F^p$; the
Plancherel--P\'olya inequalities compare the continuous and sampled
square functions, so the systems here give the same
quasi-norm. For $1<p<\infty$ the flag Littlewood--Paley theorem
gives $H_F^p=L^p$. Applying \eqref{eq:FB-embeddings} with
$s=0,q=2$ gives the first two inclusions, and
\eqref{eq:q-monotonicity} gives the last two.
\end{proof}

\section*{Funding}
No funding was received for this research.

\section*{Declaration of competing interest}
The author declares that there are no known competing financial interests or personal relationships that could have appeared to influence the work reported in this paper.

\section*{Data availability}
No datasets were generated or analysed during this study.

\section*{Declaration of generative AI and AI-assisted technologies in the manuscript preparation process}
During the preparation of this manuscript, OpenAI's ChatGPT was used to assist with language editing, manuscript organization, and the revision of mathematical arguments. The author takes full responsibility for the content of the manuscript.

\end{document}